\documentclass[12pt]{article} 
\usepackage{latexsym,amsmath}
\begin{document}

\centerline{\Large\bf Completely dissociative groupoids}\vspace{2em}

\centerline{$^1$Milton S. Braitt}

\centerline{$^2$David Hobby}

\centerline{$^2$Donald Silberger}\vspace{.5em}

\centerline{$^1$Universidade Federal de Santa Catarina (Brasil)}

\centerline{$^2$State University of New York at New Paltz (U.S.A.)}\vspace{1em}

\centerline{04 December 2010}\vspace{1em}

\centerline{\sf In memoriam Marvin E. Silberger (1898-1978), 
Ruth D. Silberger (1898-1985)}\vspace{1em}

\begin{abstract}

Consider arbitrarily parenthesized expressions on the \ $k$ \ variables \ $x_0, x_1, ..., x_{k-1}$, \ where each \ $x_i$ \ appears exactly once and in the order of their indices.  We call these expressions {\em formal $k$--products}. 
\ $F^\sigma(k)$ \ denotes the set of formal $k$--products. \ For \ $\{{\bf u},{\bf v}\}\subseteq F^\sigma(k)$, \ the claim, that \ ${\bf u}$ \ and \ ${\bf v}$ \ produce equal elements in a groupoid \ $G$ \ for all values assumed in \ $G$ \ 
by the variables \ $x_i$, \ attributes to \ $G$ \ a {\em generalized associative law}. \ Many groupoids are {\em completely dissociative}; {\it i.e.}, no generalized associative law holds for them; two examples are the groupoids on \ $\{0,1\}$ \ 
whose binary operations are implication and NAND. We prove a variety of results of that flavor.

\end{abstract}\vspace{1em}

\noindent{\large\bf \S1. \ Introduction.}\vspace{.5em}

Our preceding paper, {\bf[2]}, begins an investigation of groupoids \ ${\cal G} := \langle G;\star\rangle$ \ in which the binary operation \ $\star:G\times G\rightarrow G$ \ fails to be associative; that is, those \ ${\cal G}$ \ for which there exists an ordered triple \ $\langle g_0,g_1,g_2\rangle \in G^3$ \ with \ $(g_0\star g_1)\star g_2 \not= g_0\star(g_1\star g_2)$. \ One task that finite \ ${\cal G}$ \ of this sort inspire is to specify how many of its \ $|G|^3$ \ distinct triples do associate. Indeed, {\bf[2]} shows that, for every \ $|G|\ge 2$, \ there exists \ ${\cal G} := \langle G;\star\rangle$ \ in which every triple fails to associate.   

The failure of some triples to associate induces ambiguity in the products of longer strings as well. Thus, whereas there are only two possibly distinct products of a triple of elements in \ $G$, \ there are \ $5$ \ potentially distinct products of an ordered $4$--tuple, \ $14$ \ of a $5$--tuple, \ $42$ \ of a $6$--tuple, {\it etc.} \ The \ $(k-1)$--st  Catalan number, \ $C(k-1)$, \ happens to equal the potential diversity of the products of a $k$--tuple \ $\langle g_0,g_1,\ldots g_{k-1}\rangle \in G^k$. \ Thus a second project is suggested by the consideration of those \ ${\cal G}$ \ which are not semigroups. 

The \ $C(k-1)$ \ potentially different products of a $k$--tuple of elements in \ $G$ \ are engendered by \ $C(k-1)$ \ distinct $k$--ary formulas. We designate that formula set by \ $F^\sigma(k)$. \ Given a groupoid \ ${\cal G}$, \ some of the formulas in \ $F^\sigma(k)$ \ agree with each other everywhere on \ $G^k$, \ and so \ $F^\sigma(k)$ \ is partitioned by \ ${\cal G}$ \ into classes of mutual agreement on \ $G^k$.  \ The initiating study of such partitions is the main focus of {\bf[2]}. A detailed reminder of this work occurs in our present \S5, and it is the projected subject of a deeper study in {\bf[4]}. 

Our principal focus is upon those \ ${\cal G}$ \ for which no two distinct elements in \ $F^\sigma(k)$ \ agree everywhere on \ $G^k$, \ and for which the ``agreement classes'' in \ $F^\sigma(k)$ \ are therefore singletons. We call such groupoids $k$--dissociative.  The \ ${\cal G}$ \ that are $k$--dissociative for every \ $k$ \ we call completely dissociative.    

When \ ${\cal S} := \langle S;\,\cdot\rangle$ \ is a semigroup, there is for each positive integer \ $k$ \ and each \ $\langle s_0,s_1,\ldots,s_{k-1}\rangle \in S^k$, \ exactly one product. No ambiguity of the sort that complicates our present work arises for such products in \ ${\cal S}$, \ and so no parenthesizations are needed to distinguish one product of the $k$--tuple from another. However, failure of associativity induces ambiguity in products. So either a clutter of parentheses in the usual ``infix'' notation or a syntactic trick is needed to eliminate that ambiguity. In {\bf[2]} and here as well we opt for reverse Polish notation (rPn) whenever we deal with algebras with nonassociative binary operations. This eliminates all need for parentheses in our product expressions, except where they are used as punctuation.  

\S2 presents most of the terminology used in the main part of our article.

In \S3 we investigate $k$--dissociative and completely dissociative groupoids.  We establish a general result that aids in showing groupoids are completely dissociative, and use it to show that a surprising number of small groupoids are completely dissociative. 

In \S4 we use Birkhoff's Theorem and other tools to investigate primitive groupoids, which are minimal completely dissociative groupoids.

In \S6 we consider the construction of $k$--ary operations strings of binary operations on the same set.

Our final section, \S7, reviews and generalizes the work reported in {\bf[2]} on $k$--anti-associative groupoids, and it  introduces and poses questions about minimally $k$--associative groupoids. This subject comprises most of  {\bf[3]}.\vspace{2em}   

\noindent{\large\bf \S2. \ Our Language.}\vspace{.5em}

Henceforth \ $\omega := \{0,1,2,\ldots\}$ \ and \ ${\bf N} := \{1,2,3,\ldots\}$. \ When \ $n\in{\bf N}$ \ then \ $n$ \ also denotes the set \ $\{0,1,\ldots,n-1\}$.

For \ $\{k,n\}\subseteq{\bf N}$ \ we write \ $n^k$ \ to designate the set of all \ $k$--tuples \ of elements in \ $n$, \ and \ $n^\omega$ \ denotes the set of all infinite sequences \ $j_0j_1j_2\ldots$ \ whose terms are elements in \ $n$. \ Obviously the number of \ $k$--ary operations, \ $\phi:n^k\rightarrow n$ \ on the set \ $n$, \ is equal to the integer \ $n^{n^k}$. \ The most familiar are for \ $k=2$; \ namely, the \ $n^{n^2}$ \ distinct binary operations on the set \ $n$.
 
As in {\bf [2]}, \ we usually employ reverse Polish notation, \ rPn, \ for the \ $k$--ary operations \ $\phi:G^k\rightarrow G$. \ {\it E.g.}, for \ $k=5$ \ we write \[ \phi:\langle x_0,x_1,x_2,x_3,x_4\rangle\mapsto x_0x_1x_2x_3x_4\phi \] instead of using the more common notation \[ \phi:\langle x_0,x_1,x_2,x_3,x_4\rangle\mapsto\phi(x_0,x_1,x_2,x_3,x_4). \]

The paper {\bf [2]} discusses the set \ $F^\sigma(k)$ \ of all ``formal $k$--products", which we re-introduce in Definition 2.1, below. Our ``formal products" are special instances of what, in the more encompassing language of universal algebra, are called ``terms". An example might clarify our intent:

Consider the formal $5$--product \ ${\bf u}\,:=\,x_0x_1\bullet x_2x_3\bullet\bullet\,x_4\,\bullet\,\in\, F^\sigma(5)$.

When \ $G$ \ is a set, and if \ $\vec{\beta}$ \ is a $4$--tuple \ $\beta_0\beta_1\beta_2\beta_3$ \ with \ $\beta_i:G^2\rightarrow G$ \ for every \ $i\in 4$, \ then by our comments after Definition 2.2 below, the equalities
\[ {\bf u}^{\vec{\beta}} := x_0 x_1\bullet x_2x_3\bullet\bullet\, x_4\bullet^{\vec{\beta}} :=
x_0x_1\beta_0 x_2 x_3\beta_1\beta_2 x_4\beta_3 \] present a \ $5$--ary operation \
${\bf u}^{\vec{\beta}}:G^5\rightarrow G$ \ on \ $G$. \ The more familiar, ``infix'', notation for this $5$--product equality would be \[ {\bf u}^{\vec{\beta}}(x_0,x_1,x_2,x_3,x_4) :=
(((x_0\beta_0 x_1)\beta_2 (x_2\beta_1 x_3))\beta_3 x_4). \] The parentheses on the right of the equal sign are obligatory, since the binary operations \ $\beta_i$ \ may be highly nonassociative. However, we give the infix format a conventionally unrequired external parenthesis pair, in order to assure that the expression portrays every familiar-form product of a duple \ $\langle a,b\rangle$ \ of elements in \ $G$ \ consistently as \ $(a\beta_i b)$ \ and not as \ $a\beta_i b$, \ where \ $a\beta_i b$ \ is the colloquially ubiquitous abbreviation of \ $(a\beta_i b)$. \ This consistency, translated into rPn, simplifies discussions of subterms of formal products.

If the finite sequence \ $\vec{\beta} := \diamond\diamond\diamond\ldots\diamond$ \ of binary operations is constant, with \ $\beta_i=\diamond$ \ for all relevant \ $i$, \ then we may write \ ${\bf v}^{\vec{\beta}}$ \ more simply as \
${\bf v}^\diamond$.\vspace{.5em}

Here is a synopsis and modification of terminology introduced in {\bf[2]}:

Viewing a {\sl word} \ {\bf w} \ as a (finite or infinite) sequence, we say that \ {\bf s} \ is a {\sl subword} \ of \ 
{\bf w} \ iff \ {\bf s} \ is a subsequence of \ {\bf w}. \ A subword of \ {\bf w} \ whose letters occur consecutively 
in \ {\bf w} \ we call a {\sl segment} of \ {\bf w}; \ some authors use the word ``block'' to designate a segment in a one-letter alphabet.

An initial segment of \ {\bf w} \ we call a {\sl prefix} of \ {\bf w}; \ a terminal segment of \ {\bf w} \ we call a 
{\sl suffix} of \ {\bf w}. \ Of course a suffix of \ {\bf w} \ is infinite if and only if \ {\bf w} \ is itself
infinite.

Finite words we usually call {\sl tuples}. But from now on, infinite words will always be called {\sl sequences}, and our ``sequences''will always be infinite.

Henceforth \ $\vec{x} := x_0x_1x_2\ldots$ \ denotes a sequence of distinct variables \ $x_i$, \ and \ $\bullet$ \ denotes 
an operator symbol. The set \ $A^*$ \ is comprised of the finite words each of whose letters is an element in the infinite 
alphabet \ $A := \{\bullet,x_0,x_1,x_2,\ldots\}$. \ The objects defined in  2.1  are among the elements in the free semigroup \ $A^*$.

The binary operation on \ $A^*$ \ itself is concatenation (aka juxtaposition). However, the set \ $F^\sigma\subseteq A^*$ \ 
of formal products introduced in  2.1  is \underline{not} a subsemigroup of \ $A^*$. \ Moreover, we will introduce a binary operation \ $\odot$ \ on \ $F^\sigma$ \ which fails in the strongest way to be associative; that is, the infinite groupoid \ $\langle F^\sigma;\odot\rangle$ \ is as far from being a semigroup as possible.\vspace{.5em}

The {\sl length} of a word \ ${\bf u}\in A^*$ \ is written \ $|{\bf u}|$, \ and \ $\#({\bf u},y)$ \ denotes the number of occurrences of the letter \ $y$ \ in \ {\bf u}.\vspace{1em}

\noindent{\bf Definition 2.1.} \ Let \ ${\bf u}\in A^*$, \ and let \ $k\in{\bf N}$. \ Then \ {\bf u} \ is said to be a 
{\sl formal $k$--product} iff each of the following four criteria is satisfied:

{\bf 1.} \ \ $x_0x_1\ldots x_{k-1}$ \ is a subword of \ {\bf u}.

{\bf 2.} \ \ $|{\bf u}| = 2k-1\,$.

{\bf 3.} \ \ $\#({\bf u},\bullet)=k-1\,$.

{\bf 4.} \ \ $\#({\bf p},\bullet)<\,|{\bf p}|-\#({\bf p},\bullet)$ \
for every nonempty prefix \ {\bf p} \ of \ ${\bf u}\,$.\vspace{1em}

The expression \ $F^\sigma(k)$ \ denotes the set of all formal $k$--products, and the set of all formal products is then \
$F^\sigma := \bigcup\{F^\sigma(k): \, k\in{\bf N}\}$.

Observe that \ $F^\sigma(1) = \{x_0\}$, \ and that whenever \ $\langle{\bf a},{\bf b}\rangle\in F^\sigma(i)\times F^\sigma(j)$ 
\ then \ ${\bf ab}_i\bullet\in F^\sigma(i+j)$, \ where \ ${\bf b_i}$ \ denotes the element in \ $A^*$ \ obtained by replacing in \ {\bf b} \ the letter \ $x_t$ \ with the letter \ $x_{i+t}$ \ for each \ $t\in\omega$.

This enables us to define the binary operation \ $\odot:F^\sigma\times F^\sigma\rightarrow F^\sigma$ \ by \ ${\bf ab}\odot := {\bf ab_i}\bullet$ \ where \ ${\bf a}\in F^\sigma(i)$. \ It is routine to verify that, if \ $x_0\not={\bf w}\in F^\sigma$,
\ then there is exactly one pair \ $\langle{\bf p},{\bf s}\rangle\in F^\sigma\times F^\sigma$ \ such that \ ${\bf w} =
{\bf ps_m}\bullet\,\in F^\sigma$ \ for some \ $m\in{\bf N}$. \ That is to say, every \ ${\bf w}\in F^\sigma\setminus\{x_0\}$
\ has a unique factorization under \ $\odot$ \ into a product \ ${\bf w} = {\bf ps}\odot$ \ of two elements in \ $F^\sigma$.\vspace{1em}

\noindent{\bf Definition 2.2.} \ Let \ $\{n,k\}\subseteq{\bf N}$, \ and let \ ${\bf u}\in F^\sigma(k)$. \ For each \ $j\in k-1$ \ let \ $\beta_j:n^2\rightarrow n$, \ and let \ $\vec{\beta}$ \ be the $(k-1)$--tuple \ $\beta_0\beta_1\ldots\beta_{k-2}$. \ Then \ ${\bf u}^{\vec{\beta}}$ \ denotes the word in \ $\{\beta_0,\beta_1,\ldots\beta_{k-2},x_0,x_1,x_2,\ldots\}^*$ \
obtained by substituting the operation symbol \ $\beta_j$ \ for the \ $j$th \ occurrence of the letter \ $\bullet$ \ in the word \ ${\bf u}$, \ for each \ $j\in k-1$. \ The words \ ${\bf u}^{\vec{\beta}}$ \ are called {\sl formal} $k$--{\sl ary} $\vec{\beta}$-{\sl products}.\vspace{1em}

However, \ ${\bf u}^{\vec{\beta}}$ \ denotes also the {\sl interpretation} \ ${\bf u}^{\vec{\beta}}:G^k\rightarrow G$ \ by \ $\vec{\beta}$ \ of \ ${\bf u}$.\vspace{.5em}

\ Specifically,  \ ${\bf u}^{\vec{\beta}}$ \ denotes also the \ $k$--ary operation induced on \ $n$ \ by the formal $k$--ary
$\vec{\beta}$--product \ ${\bf u}^{\vec{\beta}}$ \ induced thus: When \ $\vec{g}:= g_0g_1\ldots\in n^\omega$ \ then \ $\vec{g}{\bf u}^{\vec{\beta}}$ \ is the element in \ $n$ \ obtained by replacing \ $x_i$ \ in \ ${\bf u}^{\vec{\beta}}$ \ with \ $g_i$ \ for each \ $i\in k$, \ and then activating the \ $k-1$ \ binary operations \ $\beta_j$.

Thus \  $F^{\sigma,\vec{\beta}}(k) := \{{\bf v}^{\vec{\beta}}:\,{\bf v}\in F^\sigma(k)\}$ \ is a set of $k$--ary operations on \ $n$. \ Obviously \ $1\le |F^{\sigma,\vec{\beta}}(k)|\le|F^\sigma(k)|$, \ and so the set \ $F^{\sigma,\vec{\beta}}$ \
is finite.\vspace{1em}

\noindent\underline{Illustrative Example 1}: \ Again let \ ${\bf u} :=
x_0x_1\bullet x_2x_3\bullet\bullet x_4\bullet\in F^\sigma(5)$. \ Let \ $n\in{\bf N}$ \ be arbitrary, let \ $\beta_j:n^2\rightarrow n$ \ be given for each \ $j\in 4$, \ and let \ $\vec{\beta} := \beta_0\beta_1\beta_2\beta_3$. \ Then \ ${\bf u}^{\vec{\beta}}$ \ denotes both the word \ $x_0x_1\beta_0x_2x_3\beta_1\beta_2x_4\beta_3$, \ and also the 
function \ ${\bf u}^{\vec{\beta}}:n^\omega\rightarrow n$ \ given by \ 
${\bf u}^{\vec{\beta}}:\vec{g}\mapsto\vec{g}{\bf u}^{\vec{\beta}} = g_0g_1\beta_0 g_2g_3\beta_1\beta_2 g_4\beta_3\in n$. \ 
Calculate the element \ $\vec{g}{\bf u}^{\vec{\beta}}\in n$ \ thus: \ Let \ $h_0:=g_0g_1\beta_0\in n$. \ Let \ $h_1:=g_2g_3\beta_1\in n$. \ Then \ $\vec{g}{\bf u}^{\vec{\beta}} = h_0h_1\beta_2g_4\beta_3 = h_2g_4\beta_3\in n$, \ where \ $h_2 := h_0h_1\beta_2$.\vspace{.5em}

When \ $\phi={\bf u}^{\vec{\beta}}$ \ then we say that \ {\bf u} \ {\em represents} \ $\phi$ \ {\it via} \ $\vec{\beta}$.\vspace{.5em}

Recall that, when \ $\beta_i  = \diamond$ \ for all \ $i$, \ where \ $\diamond$ \ is a fixed binary operation on the set \ $G$, \ we write \ ${\bf u}^{\vec{\beta}}$ \ simply as \ ${\bf u}^\diamond$. \ Incidentally, when \ $\vec{g}\in G^\omega$ \ then in the present article we write \ $\vec{g}{\bf u}^\diamond$ \ instead of \ ${\bf u}(\diamond,\vec{g})$ \ as in {\bf [2]}.

For \ $k\in {\bf N}$ \ and \ ${\bf u}\in F^\sigma(k)$ \ and \ $\vec{g}\in G^\omega$, \ note that \ 
$\vec{g}{\bf u}^{\vec{\beta}}$ \ is determined by the length--$k$ prefix of \ $\vec{g}$.  \ The ``extra" terms in \ $\vec{g}$ \ simplify our notation.

Two additional conventions: When \ $y\in G$ \ and \ $\vec{g}:=g_0g_1g_2\ldots\in G^\omega$ \ is the sequence such that \ $g_t = y$ \ for all \ $t\in\omega$ \ then we may write \ $\vec{g}$ \ instead as \ $\vec{y}$, \ where \ $\vec{y} := yyy\ldots$ \ That is to say, when \ $y\in G$ \  then \ $\vec{y} := yyy\ldots\in G^\omega$.

For \ $\langle G;\diamond\rangle$ \ a groupoid, when \ $\vec{g} = g_0g_1g_1\ldots\in G^\omega$ \ 
and \ $m\in\omega$ \ then \ $\vec{g}_m$ \ denotes the infinite suffix \ $g_mg_{m+1}g_{m+2}\ldots$ \ 
of \ $\vec{g}$. \ Thus \ $\vec{g}{\bf u_m}^\diamond = \vec{g}_m{\bf u}^\diamond$.\vspace{2em}

\noindent{\large\bf\S3. \ Complete dissociativity.}\vspace{.5em}

For \ ${\cal G} := \langle G;\diamond\rangle$ \ a groupoid, {\bf[2]} defines an equivalence relation \ 
$\approx_{\cal G}$ \ on the set \ $F^\sigma(k)$ \ by \ ${\bf u}\approx_{\cal G}{\bf v}$ \ iff \
$\vec{g}{\bf u}^\diamond=\vec{g}{\bf v}^\diamond$ \ for each \ $\vec{g}\in G^\omega$; \ that is, iff \ 
${\bf u}^\diamond = {\bf v}^\diamond$. \ \ $F^\sigma(k)/{\cal G}$ \ denotes the family of all \ 
$\approx_{\cal G}$ \ equivalence classes \ $[{\bf u}]_{\cal G}$.\vspace{.5em}

If \ $\vec{g}{\bf u}^\diamond \not= \vec{g}{\bf v}^\diamond$ \ then we say that \ $\vec{g}$ \ {\sl separates} \ 
${\bf u}^\diamond$ \ from \ ${\bf v}^\diamond$.\vspace{1em}

\noindent{\bf Definition.} \ \ $\langle G;\diamond\rangle$ \ is said to be $k$--{\sl dissociative} iff 
every pair \ ${\bf u}\not={\bf v}$ \ of elements in \ $F^\sigma(k)$ \ can be separated by a sequence in 
\ $G^\omega$. \ We call \ $\langle G;\diamond\rangle$ \ {\sl completely dissociative} iff \ 
$\langle G;\diamond\rangle$ is \ $k$--dissociative for every \ $k\ge 3$.\vspace{1em}

The $3$--dissociative groupoids are those which are not semigroups.\vspace{.5em}

Some workers call \ ${\bf u}^\diamond = {\bf v}^\diamond$ \ an {\sl identity} in \ ${\cal G}$. \ Our formal products 
also correspond to specific kinds of {\sl terms}; namely, those where the variables \ $x_0,x_1,x_2,\ldots x_{k-1}$ \ 
each appear exactly once, and in the order given here. We will hereafter call terms of this sort (interpreted) formal 
products.

\ Accordingly, \ $\langle G;\diamond\rangle$ \ is $k$--dissociative if and only if no nontrivial $k$--ary 
identity is satisfied between formal products interpreted in \ ${\cal G} := \langle G;\diamond\rangle$.\vspace{.5em}

\noindent{\bf N.B.} \ {\sf Our completely dissociative \ $\langle G;\diamond\rangle$ \
were called ``completely free" in } {\bf[1]}.\vspace{.5em}

We begin with a few remarks. For every integer \ $k\ge 3$, \ notice that:\vspace{.5em}

${\cal G} := \langle G;\diamond\rangle$ \ is $k$--dissociative if, and only if, \ 
$|F^\sigma(k)/{\cal G}|=|F^\sigma(k)|$ \ and each \ $\approx_{\cal G}$ \ equivalence class is a singleton \ 
$[{\bf v}]_{\cal G} = \{{\bf v}\}\subseteq F^\sigma(k)$.

Both isomorphism and anti-isomorphism respect $k$--dissociativity.

${\cal G}:=\langle G;\diamond\rangle$ \ is $k$--dissociative if \ ${\cal G}$ \ has a $k$--dissociative 
subgroupoid.

If a component groupoid of a cartesian product groupoid is $k$--dissociative then the product groupoid too 
is $k$--dissociative.

If \ ${\cal G}$ \ has a $k$--dissociative homomorphic image, \ ${\cal G}$ \ is $k$--dissociative.\vspace{1em}
 
While the five statements above can all be verified directly, they are also consequences of Birkhoff's Theorem 
in Universal Algebra.  This will be discussed in more detail in \S4.\vspace{.5em}

\noindent{\bf Question.} \ For each \ $k\ge 3$ \ is there a $k$--dissociative groupoid which is not $(k+1)$--dissociative?\vspace{1em}

\noindent{\bf Theorem 3.1.} \ {\sl If \ $\langle G;\diamond\rangle$ \ is $k$--dissociative, then \ 
$\langle G;\diamond\rangle$ \ is $j$--dissociative for all \ $j\in\{3,4,\ldots,k-1\}$. }\vspace{1em}

\noindent{\bf Proof.} \ We apply induction to the contrapositive. The claim holds for \ $k=3$. \ Pick \ $k\ge 3$. 
\  Suppose that \ ${\bf a}^\diamond = {\bf b}^\diamond$ \ for some formal $k$--products \ ${\bf a} \not= {\bf b}$. 
\ Then of course \ $\{{\bf a}x_k\bullet,\,{\bf b}x_k\bullet\}\subseteq F^\sigma(k+1)$ \ and \ 
${\bf a}x_k\bullet \not= {\bf b}x_k\bullet\,$. \ Let \ $\vec{g}\in G^\omega$ \ be arbitrary. Since \ 
$\vec{g}{\bf a}^\diamond = \vec{g}\,{\bf b}^\diamond$, \ we have \ $\vec{g}{\bf a}x_k\bullet^\diamond =
\vec{g}{\bf a}^\diamond g_k\diamond = \vec{g}\,{\bf b}^\diamond g_k\diamond = \vec{g}\,{\bf b}x_k\bullet^\diamond\,$. 
\ Thus \ ${\bf a}x_k\bullet^\diamond = {\bf b}x_k\bullet^\diamond$.\hfill$\rule{.5em}{.8em}$\vspace{1em}

Conjecture {\bf [2}:3.13{\bf]} fails. Of the \ $16$ \ binary operation tables on the set \ $2$, \ eight are of 
semigroups. We call the tables themselves ``concrete" semigroups.

Of the eight concrete nonsemigroups on the set \ $2$, \ our computer verifies that two fail to be $4$--dissociative. 
Theorem 3.6 establishes that the remaining six are completely dissociative.\vspace{.5em}

In order easily to keep track of the entities we intend to discuss, we use a natural nomenclature denoting specific 
binary operation tables on the sets \ $n\in{\bf N}$. \ This dictionary of binary operation tables arises from the 
fact that there are \ $n^{n^2}$ \ such tables for the set \ $n$. \ If practicality is ignored, we can express
each such table as a Hindu-Arabic numeral to the base \ $n$ \ of an integer \ $j\in n^{n^2}$. \ This integer \ $j$ 
\ names and specifies the binary operation table \ $n_j$, \ provided that one reads the table as you are reading 
the page before you now, from upper left to lower right. The numeral must be \ $n^2$ \ digits long, and so we may 
need a prefix of \ $0$s in order to assure a word that is exactly \ $n^2$ \ base--$n$ \ digits long.

$\underline{\cal G}(2) := \{2_0,2_1,2_2,\ldots,2_{15}\}$ \ is the family of all \ $2^{2^2} = 16$ \ distinct such 
tables on the universe \ $2:=\{0,1\}$. \ The groupoid \ $2_j$ \ is a semigroup if and only if \ 
$j\in\{0,1,3,5,6,7,9,15\}$. \ Our computer verifies that neither \ $2_{10}$ \ nor \ $2_{12}$ \ is $4$--dissociative,
although not one of the eight triples \ $\langle a,b,c\rangle\in 2^3$ \ associates in \ $2_{10}$ \ or in \ $2_{12}$.\vspace{.5em}

Initially we established by induction the complete dissociativity of several groupoids \ ${\cal G}$. \ Then, 
similarities in those proofs suggested a comprehensive result, Theorem 3.2, below, which requires additional 
terminology.

When a $j$--tuple \ $\vec{r}\in G^j$ \ occurs \ $m$ \ times consecutively in a sequence \ $\vec{g}\in G^\omega$, \ 
we write the resulting  $mj$--tuple segment of \ $\vec{g}$ \ as \ ${\vec{r}\,}^m$. \ This generalizes \ $y^m$, \ 
which denotes the \ $m$--tuple \ $yy\ldots y\in G^m$ \ when \ $y\in G$.\vspace{1em}

\noindent{\bf Definitions.} \ For \ $\vec{g} \in G^\omega$ \ where \ ${\cal G} := \langle G;\diamond\rangle$ \ 
is a groupoid, and for \ $S\subseteq G$, \ we say that \ $\vec{g}$ \ {\em yields} \ $S$ \ iff \ 
$\vec{g}{\bf u}^\diamond \in S$ \ for every \ ${\bf u} \in F^\sigma$.

Given \ $U \subseteq F^\sigma$, \ we say \ $\vec{g}$ \ yields \ $S$ \ {\em on} \ $U$ \ iff \ $\vec{g}{\bf u}^\diamond \in S$ \ for all \ ${\bf u}\in U$. 

The set \ $S$ \ is called {\em yieldable} iff there is some \ $\vec{g}$ \ which yields \ $S$. \ If \ $\vec{g}$ \ 
yields \ $\{a\}$ \ then we say that \ $\vec{g}$ \ yields \ $a$ \ and that \ $a$ \ is yieldable.  

For \ $i\in{\bf N}$ \ and \ ${\bf u}\in F^\sigma$, \ we call \ ${\bf u}$ \ an {\em $i$--split} iff its 
unique factorization in \ $\langle F^\sigma;\odot\rangle$ \ is \ ${\bf u} = {\bf ab}\odot$ \ with \ 
${\bf a}\in F^\sigma(i)$.\vspace{1em}

\noindent{\bf Theorem 3.2.} \ {\sl Let \ ${\cal G} := \langle G;\diamond\rangle$ \ be a groupoid, let \ 
$T \subseteq G$ \ with \ $|T|\ge 2$, \ and let the following three conditions hold.\vspace{.2em}

{\bf 1.} {\sf Left Separation:} \ If \ $\{x,y\}\subseteq T$ \ with \ $x\not= y$ \ then there is a 
yieldable \ $L_{x,y}\subseteq G$ \ such that \ $s\diamond x\not=s'\diamond y$ \ with \ 
$\{s\diamond x,s'\diamond y\}\subseteq T$ \ for every \ $\langle s,s'\rangle \in L_{x,y}\times L_{x,y}\,$.

{\bf 2.} {\sf Right Separation:} \ If \ $\{x,y\}\subseteq T$ \ with \ $x\not= y$ \ then there is a 
yieldable \ $R_{x,y}\subseteq G$ \ such that \ $x\diamond s\not= y\diamond s'$ \ with \ 
$\{x\diamond s,y\diamond s'\}\subseteq T$ \ for every \ $\langle s,s'\rangle \in R_{x,y}\times R_{x,y}$.

{\bf 3.} {\sf Split Separation:} \ For all \ $\{i,j\}\in{\bf N}$ \ with \ $i\not=j$, \ there exist nonempty 
disjoint subsets \ $A$ \ and \ $B$ \ of \ $T$ \ and some \ $\vec{g}\in G^\omega$ \ which yields \ $A$ \ on 
the set of all $i$--splits and which yields \ $B$ \ on the set of all $j$--splits.\vspace{.3em}

Then \ ${\cal G}$ \ is completely dissociative.}\vspace{1em}

\noindent{\bf Proof.} \ We will prove by induction that \ ${\cal G}$ \ is \ $k$--dissociative for every \ 
$k\in{\bf N}$. \ This claim holds trivially for \ $k\in\{1,2\}$. \ Note that the $1$--split \ 
$x_0x_1x_2\bullet\bullet$ \ and the $2$--split \ $x_0x_1\bullet x_2\bullet$ \ are the only two elements in \ 
$F^\sigma(3)$. \ By 3.2.3, there exist \ $A\subseteq T$ \ and \ $B\subseteq T$ \ and \ $\vec{g}\in G^\omega$ \ 
with \ $A\cap B = \emptyset$, \ and such that \ $a = \vec{g}x_0x_1x_2\bullet\bullet^\diamond$ \ for some \ 
$a\in A$ \ and \ $b = \vec{g}x_0x_1\bullet x_2\bullet^\diamond$ \ for some \ $b\in B$. \ So \  ${\cal G}$ \ is 
$k$--dissociative if \ $k\in\{1,2,3\}$. \ The basis is done.

Now choose \ $k\in\{3,4,5,\ldots\}$, \ and suppose that \ ${\cal G}$ \ is \ $v$--dissociative for every \ 
$v\in\{1,2,\ldots,k\}$. \ Pick any \ $\{{\bf w},{\bf w'}\}\subseteq F^\sigma(k+1)$ \ with \ ${\bf w} \not= 
{\bf w'}$. \ These formal $(k+1)$--products have unique factorizations \ ${\bf w} = {\bf ab}\odot$ \ and \ 
${\bf w'} = {\bf a'b'}\odot$ \ in the groupoid \ $\langle F^\sigma;\odot\rangle$, \ where \ 
$\langle{\bf a},{\bf a'}\rangle \in F^\sigma(i)\times F^\sigma(i')$ \ for some \ 
$\{i,i'\}\subseteq\{1,2,\ldots,k\}$. \ Without loss of generality, take it that \ $i\le i'$.

If \ $i<i'$ \ then by 3.2.3 there exist disjoint subsets \ $A$ \ and \ $B$ \ of \ $T$, \ and a sequence \ 
$\vec{g}\in G^\omega$, \ such that \ $\vec{g}{\bf w}^\diamond \in A$ \ while \ $\vec{g}{\bf w'}^\diamond \in B$, 
\ whence \ $\vec{g}{\bf w}^\diamond \not= \vec{g}{\bf w'}^\diamond$, \ and so \ ${\bf w}^\diamond \not= 
{\bf w'}^\diamond$. \ So we may take it that \ $i=i'$. \ Since \ ${\bf w} \not= {\bf w'}$, \ either \ 
${\bf a}\not={\bf a'}$ \ or \ ${\bf b} \not= {\bf b'}$.\vspace{.5em}

\underline{Case}: \ ${\bf a} \not= {\bf a'}$. \ By the inductive hypothesis, there exists \ $\vec{g}\in G^\omega$ 
\ with \ $\vec{g}{\bf a}^\diamond \not= \vec{g}{\bf a'}^\diamond$ \ and \ 
$\{\vec{g}{\bf a}^\diamond,\vec{g}{\bf a'}^\diamond\}\subseteq T$. \ Let \ $x := \vec{g}{\bf a}^\diamond$ \ and \ 
$x' := \vec{g}{\bf a'}^\diamond$. \ By 3.4.2 there exists \ $\vec{h}\in G^\omega$ \ such that \ $xs\diamond$ \ and \ $x's'\diamond$ \ are distinct elements in \ $T$, \ where \ $s := \vec{h}{\bf b}^\diamond$ \ and \ $s' := \vec{h}{\bf b'}^\diamond$. \ Clearly we may suppose \ $\vec{g_i} = \vec{h}$. \ Thus \ $\vec{g}{\bf w}^\diamond = \vec{g}{\bf ab}\odot^\diamond = \vec{g}{\bf ab_i}\bullet^\diamond = \vec{g}{\bf a}^\diamond\vec{g}{\bf b_i}^\diamond\diamond = 
\vec{g}{\bf a}^\diamond\vec{g_i}{\bf b}^\diamond\diamond = xs\diamond \not= x's'\diamond = \vec{g}{\bf a'}^\diamond\vec{g_i}{\bf b'}^\diamond\diamond = \vec{g}{\bf a'}^\diamond\vec{g}{\bf b'_i}^\diamond\diamond = \vec{g}{\bf a'b'_i}\bullet^\diamond = \vec{g}{\bf a'b'}\odot^\diamond = \vec{g}{\bf w'}^\diamond$, \ and so \  
${\bf w}^\diamond \not= {\bf w'}^\diamond$.\vspace{.5em}

\underline{Case}: \ ${\bf b} \not= {\bf b'}$. \ By the inductive hypothesis, there exists \ $\vec{r}\in G^\omega$ 
\ with \ $\vec{r}\,{\bf b}^\diamond \not= \vec{r}\,{\bf b'}^\diamond$ \ and \ 
$\{\vec{r}\,{\bf b}^\diamond,\vec{r}\,{\bf b'}^\diamond\}\subseteq T$. \ Let \ $x := \vec{r}\,{\bf b}^\diamond$ \ and \ 
$x' := \vec{r}\,{\bf b'}^\diamond$. \ By 3.4.1 there exists \ $\vec{g}\in G^\omega$ \ such that \ $sx\diamond$ \ 
and \ $s'x'\diamond$ \ are distinct elements in \ $T$, \ where \ $s := \vec{g}{\bf a}^\diamond$ \ and \ $s' := 
\vec{g}{\bf a'}^\diamond$. \ We may suppose that \ $\vec{g_i} = \vec{r}$. \ Thus \ $\vec{g}{\bf w}^\diamond = 
\vec{g}{\bf ab}\odot^\diamond = \vec{g}{\bf ab_i}\bullet^\diamond = 
\vec{g}{\bf a}^\diamond\vec{g}{\bf b_i}^\diamond\diamond = 
\vec{g}{\bf a}^\diamond\vec{g_i}{\bf b}^\diamond\diamond = sx\diamond \not= s'x'\diamond = 
\vec{g}{\bf a'}^\diamond\vec{g_i}{\bf b'}^\diamond\diamond = 
\vec{g}{\bf a'}^\diamond\vec{g}{\bf b'_i}^\diamond\diamond = \vec{g}{\bf a'b'_i}\bullet^\diamond = 
\vec{g}{\bf a'b'}\odot^\diamond = \vec{g}{\bf w'}^\diamond$, \ and so \ 
${\bf w}^\diamond \not= {\bf w'}^\diamond$.\vspace{.5em}    

Thus \ ${\cal G}$ \ is $(k+1)$--dissociative. So \ ${\cal G}$ \ is completely  dissociative.\hfill$\rule{.5em}{.8em}$\vspace{1em} 

Although Theorem 3.2 is quite general, it is often used in a simple way. For instance, if \ $a\in G$ \ 
is idempotent then surely \ $a$ \ is yieldable, for we can let \ $\vec{g} = a^\omega$. \ Likewise, if \ 
$x\mapsto xa\diamond$ \ is a permutation of \ $T$ \ then Right Separation is shown by setting \ $R_{x,y} 
:= \{a\}$ \ for all \ $x\not= y$. \ And, if \ $\diamond$ \ is commutative that Right Separation is 
equivalent to Left Separation. It may happen also that \ $T = G$ \ instead of \ $T\subset G$.

Some situations arise repeatedly when we argue that \ $\vec{g}$ \ yields a particular set.  Note
that if \ $H$ \ is a subgroupoid of \ $\cal{G}$, \ and if \ $g_0g_1\ldots g_{k-1} \in H^k$ \ for some \ 
$k\in{\bf N}$, \ then \ $\vec{g}$ \ yields \ $H$. \ Another situation arises when \ $\{a,b\}$ \ is a 
$2$--element subgroupoid of \ $\cal{G}$ \ with \ $aa\diamond = a$ \ and with \ $ab\diamond = ba\diamond 
= bb\diamond = b$; \ that is, the set \ $\{a,b\}$ \ {\em forms a semilattice}, and when also \ 
$a^\omega \not= \vec{g}\in\{a,b\}^\omega$, \ then \ $\vec{g}$ \ yields \ $b$, \ because \ $b$ \ is an 
absorptive element in subgroupoid \ $\{a,b\}$.  

As our first example of the use of Theorem 3.2, we prove the following.\vspace{1em}

\noindent{\bf Theorem 3.3.} \ {\sl The groupoid \ ${\cal B} := \langle 4;\beta\rangle$, \ below, is completely
dissociative.}\vspace{1em} 

\noindent{\bf Proof.} \ Let \ ${\cal B} := \langle 4;\beta\rangle$ \ be the groupoid whose binary operation \ 
$\beta$ \ is given by the following table:\vspace{1em}

\setlength{\unitlength}{1em}

\begin{picture}(20,11)

\thicklines 

\put(9.7,9){\line(1,0){9}} 

\put(11,2.2){\line(0,1){8}} 

\put(9.8,9.4){$\beta$}

\put(12,9.3){\bf 0} 

\put(14,9.3){\bf 1} 

\put(16,9.3){\bf 2} 

\put(18,9.3){\bf 3} 

\put(10,7.5){\bf 0} 

\put(10,5.8){\bf 1} 

\put(10,4.1){\bf 2} 

\put(10,2.4){\bf 3} 

\put(12,7.5){$0$} 

\put(14,7.5){$1$} 

\put(16,7.5){$2$} 

\put(18,7.5){$3$} 

\put(12,5.8){$1$} 

\put(12,4.1){$2$} 

\put(12,2.4){$3$} 

\put(14,5.8){$1$} 

\put(16,5.8){$3$} 

\put(18,5.8){$2$} 

\put(14,4.2){$3$} 

\put(16,4.2){$2$} 

\put(18,4.2){$1$} 

\put(14,2.4){$2$} 

\put(16,2.4){$1$} 

\put(18,2.4){$3$} 

\put(14.5,0){\large${\cal B}$}

\end{picture}

\vspace{1em}

To apply 3.2, let \ $T := 4$. \ Note that \ $0$ \ is yieldable since \ $0 = 00\beta$, \ using rPn  
language. Left and Right Separabilities are equivalent, since \ ${\cal B}$ \ is abelian. \ \ $0$ \ 
is an identity element; for each \ $x\not=y$, \ let \ $L_{x,y} = R_{x,y} = \{0\}$.

To show Split Separability, let \ $1\le i < j < k$, \ and choose \ $\vec{g} := 
0^{i-1}120^{j-i-1}30^{k-j-1}\vec{g_k}$. \ Let \ $A := \{1\}$ \ and \ $B := \{3\}$. \ For an  
$i$--split \ ${\bf w}={\bf ab}\odot\in F^\sigma(k)$ \ with \ $\langle{\bf a},{\bf b}\rangle\in F^\sigma(i)\times F^\sigma(k-i)$, \ we compute that \ $\vec{g}{\bf w}^\beta = \vec{g}{\bf ab}\odot^\beta =  
\vec{g}{\bf ab_i}\bullet^\beta = \vec{g}{\bf a}^\beta\vec{g}{\bf b_i}^\beta\beta = 
\vec{g}{\bf a}^\beta\vec{g_i}{\bf b}^\beta\beta = (0^{i-1}1\vec{g_i}){\bf a}^\beta\vec{g_i}{\bf b}^\beta\beta 
= 1(\vec{g_i}{\bf b}^\beta)\beta = 1(20^{j-i-1}30^{k-j-1}\vec{g_k}{\bf b}^\beta)\beta = 123\beta\beta = 
11\beta = 1\in A$. \ That is to say, \ $\vec{g}$ \ yields \ $A$ \ on the set of all $i$--splits in \ $F^\sigma(k)$. 
\ Similarly \ $\vec{g}$ \ yields \ $12\beta3\beta = 3\in B$ \ on the set of all \ $j$--splits in \ $F\sigma(k)$. \ So if \ ${\bf w'}\in F^\sigma(k)$ \ is a $j$--split then \ $\vec{g}{\bf w}^\diamond \not= \vec{g}{\bf w'}^\diamond$. \ That is, \ $\vec{g}$ \ separates \ ${\bf w}$ \ and  \ ${\bf w'}$.\hfill$\rule{.5em}{.8em}$\vspace{1em}

Every finite nonempty semigroup contains an idempotent; all elements in \ ${\cal B}$ \ are idempotents. The 
commutativity of \ $\beta$ \ forces association on each of the sixteen triples \ $\langle x,y,x\rangle \in 4^3$ 
\ out of a total of sixty four triples. It is as if \ ${\cal B}$ \ were trying to be a semigroup, and it is 
surprising that such a groupoid should be completely dissociative.

${\cal B}$ is interesting also because of its subgroupoid \ $\{1,2,3\}$, \ isomorphic to the \ ${\cal E}$ \ that  
we discuss below. \ ${\cal E}$ \ falls just short of complete dissociativity.  It is curious that removing  
the identity element, \ $0$, \  from \ ${\cal B}$ \ destroys complete dissociativity.\vspace{.5em}

We need a few observations about the concrete groupoid \ $2_{13} := \langle 2;\diamond\rangle$, \ which earns 
from us the name of ``implication''. It is given by the base--$2$ Hindu-Arabic numeral \ $1101$; \ that is to say, \ 
$\diamond$ \ is defined by \ $00\diamond = 01\diamond = 11\diamond = 1$ \ but \ $10\diamond = 0$. \ A standard 
interpretation of \ $0$ \ as ``false'' and \ $1$ \ as ``true'', would justify our interpreting \ $\diamond$ \ as \ 
$\Rightarrow$, \ presented here in reverse Polish notation.\vspace{1em}

\noindent{\bf Lemma 3.4.} \ {\sl Let \ $\vec{g}\in 2^\omega$, \ let \ $k\in{\bf N}$, \ and let \ 
${\bf u}\in F^\sigma(k)$. \ Then the groupoid \ $2_{13}$ \ satisfies the following conditions.

{\bf 1.} \ If \ $k\ge1$ \ and \ $\vec{g} = 1^k\vec{g}_k$ \ then \ $\vec{g}{\bf u}^\diamond = 1$.

{\bf 2.} \ If \ $k \ge 1$ \ and \ $\vec{g} = 1^{k-1}0\vec{g}_k$ \ then \ $\vec{g}{\bf u}^\diamond = 0$.

{\bf 3.} \ If \ $k \ge 2$ \ and \ $j \le k-2$ \ and \ $\vec{g} = 1^j01^{k-j-1}\vec{g}_k$ \ then  
\ $\vec{g}{\bf u}^\diamond = 1$. }\vspace{1em}

\noindent{\bf Proof.} \ The claim 3.4.1  follows from the fact that \ $1$ \ is an idempotent.

We prove 3.4.2 by induction.  The basis step, where \ $k=1$, \ is obvious. So pick \ $k\ge1$, \ and suppose that 
the lemma holds for all \ $i\in\{1,2,\ldots,k\}$. \ Let \ ${\bf v}={\bf ab}\odot\in F^\sigma(k+1)$, \ with \ 
$\langle{\bf a},{\bf b}\rangle\in F^\sigma(i)\times F^\sigma(k+1-i)$ \ for some \ $i\in\{1,2,\ldots,k\}$. 
\ Let \ $\vec{g} := 1^k0\vec{g}_{k+1}\in 2^\omega$. \ Since \ $1\le k+1-i\le k$, \ and since \ $\vec{g}_i=1^{k-i}0\vec{g}_{k+1}$, \ we have by the inductive hypothesis that \ $\vec{g}_i{\bf b}^\diamond = 0$. 
\ Thus \ $\vec{g}{\bf v}^\diamond = \vec{g}{\bf ab}\odot^\diamond = 
\vec{g}{\bf a}^\diamond\,\vec{g}\,{\bf b}_i^\diamond\diamond =
\vec{1}{\bf a}^\diamond\,\vec{g}_i{\bf b}^\diamond\diamond = 1\,0\,\diamond = 0$. \ So 3.4.2 follows.

The basis of an inductive proof of 3.4.3 involves \ $k = 2$ \ and \ $j = 0$ \ and \ $\vec{g} = 01\vec{g}_2\in 2^\omega$. 
\ For the only \ ${\bf u}\in F^\sigma(2)$ \ we then get \ $\vec{g}{\bf u}^\diamond = 01\diamond = 1$.  

Now pick \ $k\ge 2$, \ and let \ ${\bf v}={\bf ab}\odot\in F^\sigma(k+1)$ \ be a $j$--split with \ $j\le k-2$. \ 
Let \ $\vec{g} = 1^j01^{k-j}\vec{g}_{k+1}$. \ Suppose, for all \ $t\in\{1,2,\ldots,k\}$ \ and \ $i\le t-2$, \ 
that \ $1^i01^{t-i-1}\vec{g}_t$ \ yields \ $1$ \ on the set of all $i$--splits in \ $F^\sigma(t)$. \ 
Then \ $\vec{g}_j{\bf b}^\diamond = 1$, \ either by the inductive hypothesis, or by 3.4.1.  
In any event, \ $\vec{g}{\bf v}^\diamond \in \{01\diamond,11\diamond\} = \{1\}$.\hfill$\rule{.5em}{.8em}$\vspace{1em}

\noindent{\bf Theorem 3.5.} \ {\sl The implication groupoid \ $2_{13}$ \ is completely dissociative.}\vspace{1em}

\noindent{\bf Proof.} \ Take \ $T := \{0,1\}$, \ and use 3.2.  For Left Separation, let \ $L_{x,y} = \{1\}$ \ for all 
\ $x\not= y$; \ this set is yieldable by Lemma 3.4.1. For Right Separation, let \ $R_{x,y} = \{0\}$ \ for all \ 
$x\not= y$; \ this set is yieldable by 3.4.2.

To show Split Separation, let \ $1\le i<j\le k$, \ let \ $\vec{g} := 1^{i-1}01^{k-i}0\vec{g}_{k+1}$, \ let \ 
$\bf{ab\odot}\in F^\sigma(k+1)$ \ be an $i$--split, and let \ ${\bf a'b'}\odot\in F^\sigma(k+1)$ \ be a $j$--split. 
Then by 3.4 we get that \ $\vec{g}{\bf ab\odot}^\diamond = 
(1^{i-1}0\vec{g}_i{\bf a}^\diamond)(1^{k-i}0\vec{g}_{k+1}{\bf b}^\diamond)\diamond =
00\diamond = 1$, \ while on the other hand \ $\vec{g}{\bf a'b'\odot}^\diamond = 
(1^{i-1}01^{j-i}{\bf a}^\diamond)( 1^{k-j}0{\bf b}^\diamond)\diamond = 10\diamond = 0$.\hfill$\rule{.5em}{.8em}$\vspace{1em}

To complete our determination of the completely dissociative \ $2_j$, \ we deal in 3.6 with the ``NAND'' groupoid, \ $2_{14}$. \ Our proof Theorem 3.6 is aberrant, in that it does not use Theorem 3.2, since Theorem 3.2 requires some \ $\vec{g}$ \ that reliably yield particular sets, and this seems not feasible with the idempotent-free groupoid \ $2_{14}$.\vspace{1em}

\noindent{\bf Theorem 3.6.} \ {\sl The NAND groupoid \ $2_{14}$ \ is completely dissociative.}\vspace{1em}

\noindent{\bf Proof.} \ For reference, the table of $2_{14}$ is 
given below.  

\setlength{\unitlength}{1em}

\begin{picture}(20,8)

\thicklines 

\put(10.7,6){\line(1,0){5}} 

\put(12,2.6){\line(0,1){4.5}} 

\put(11,6.4){$\star$}

\put(13,6.3){\bf 0} 

\put(15,6.3){\bf 1} 

\put(11,4.5){\bf 0} 

\put(11,2.8){\bf 1} 

\put(13,4.5){$1$} 

\put(15,4.5){$1$} 

\put(13,2.8){$1$} 

\put(15,2.8){$0$} 

\put(13,1.3){\large $2_{14}$}
 
\end{picture}

The binary operation \ $\star$ \ of \ $2_{14}$ \ is equivalent to an expression in the standard boolean 
algebra on \ $2 := \{0,1\}$. \ The binary operations of this boolean algebra are {\em join} or {\em sum}, 
written \ $\vee$, \ and {\em meet} or {\em product}, written \ $\wedge$, \ and its unary operation is 
{\em complement}, written \ $\prime\,$. \ With this notation, we have that \ $xy\star = x'\vee y'$ \ for 
all \ $\langle x,y\rangle\in 2\times 2$, \ read as  \ ``$x$ NAND $y$''. 

Our proof will proceed {\it via} boolean algebra expressions that are equivalent to formal products.  
These expressions will be reduced to a standard form similar to disjunctive normal form.  The following 
terminology is due mainly to W. V. Quine; {\it viz} {\bf [8]} or Chapter XIV of {\bf [9]}. However, our 
presentation will be self-contained. 

Expressions will be built up out of variables; a {\em literal} will be either a single variable \ $x_i$ \ 
or its complement \ $x'_i$. \ A {\em fundamental formula} is either a single literal or a conjunction 
of literals with no repeated variables. A formula \ $\Phi$ \ is {\em normal} if it is either fundamental or 
a disjunction of fundamental formulas.  In the latter case, the fundamental formulas are {\em clauses} of $\Phi$.

A formula \ $\Theta$ \ is said to {\em imply} a formula \ $\Phi$ \ iff every uniform assignment of values 
to the variables in the formulas makes \ $\Phi$ \ equal to \ $1$ \ if it makes \ $\Theta$ \ equal to \ $1$; \ 
we then call \ $\Theta$ \ an {\em implicant} of \ $\Phi$. \ A {\em prime implicant} of \ $\Phi$ \ is a fundamental 
formula that implies \ $\Phi$, \ but fails to do so if any of its literals is removed.

The formulas are also called ``Sum of Product" or SoP forms.  Our focus is upon a special SoP form, called the 
{\em complete sum} form. The complete sum of a formula \ $\Phi$, \ which is equivalent neither to \ $0$ \ nor to \ $1$, 
\ is defined to be the disjunction of all its prime implicants. (A formula equivalent to \ $0$ \ has no implicants; a 
formula equivalent to \ $1$ \ has an ``empty product" as its sole prime implicant.  We avoid these {\em trivial} cases.)  
It is easy to recognize the implicants of a nontrivial formula, and the prime implicants are clearly identifiable.  The 
complete sum form of a nontrivial formula is unique, up to the order of clauses and of literals within clauses.

For example, in the formula \ $\Phi := (x\wedge y)\vee(x\wedge y')\vee z\vee(x'\wedge y\wedge z)$, \ each of its four 
clauses \ $x\wedge y$, \ $x\wedge y'$, \ $z$ \ and \ $x'\wedge y\wedge z$ \ are implicants of \ $\Phi$, \ as are such 
fundamental formulas such as \ $x\wedge z$ \ and \ $x'\wedge y'\wedge z$. \ The clause \ $z$ \ is a prime implicant of \ 
$\Phi$, \ but \ $x\wedge y$ \ and \ $x\wedge y'$ \ are not -- they can be {\em combined} into the fundamental formula \ 
$x$. \ The final clause \ $x'\wedge y\wedge z$ \ also fails to be a prime implicant; it is {\em subsumed} by \ $z$, \ 
and thus can be {\em deleted}. So the prime implicants of \ $\Phi$ \ are \ $x$ \ and \ $z$, \ and the complete sum form 
of \ $\Phi$ \ is \ $x\vee z$. \ Quine attributes this process of combining and deleting clauses to Samson and Mills, and presents a proof that it always yields our complete sum form of a nontrivial formula. (Quine calls our complete sum form 
of a formula ``the alternation of its prime implicants".)  His proof is sometimes called {\em Quine's Theorem}; it states 
that a formula is in complete sum form if and only if no clauses can be combined or deleted.  It could be used to simplify 
the proofs of Claims 2 and 3, below.\vspace{.5 em}

\noindent\underbar{Claim 1}: \ If \ ${\bf u}\in F^\sigma(k)$ \ then there exists \ $\vec{g}\in 2^\omega$ \ such that \ 
$\vec{g}{\bf u}^\star = 0$ \ and \ $\vec{r}\in 2^\omega$ \ such that \ $\vec{r}{\bf u}^\star = 1$. \ These 
evaluations depend on all \ $x_i$ \ for \ $i\in k$.\vspace{.5 em}

The claim is obvious for \ $k=1$. \ If it holds for \ ${\bf u},{\bf v}\}\subseteq F^\sigma$, \ it holds for \ 
${\bf uv}\odot$. \ So induction establishes Claim 1, none of our \ ${\bf w}\in F^\sigma$ \ are trivial, and we can  
restrict our focus to the complete sum form of \ ${\bf w}$.\vspace{.5em}

\noindent \underbar{Claim 2}: \ If \ $p = st\star\in F^{\sigma,\star}$ \ then the complete sum form of \ $p$ \ is 
equal to the join of the complete sum forms of \ $s'$ \ and \ $t'$.\vspace{.5em}

To prove this, first observe that \ $s$, \ $t$ \ and \ $p$ \ are nontrivial by Claim 1. So \ $s'$ \ and \ $t'$ \ also 
are nontrivial. So the complete sum forms of \ $p$, \ of \ $s$, \ of \ $s'$, \ of \ $t$, \ and of \ $t'$ \ all exist. 
Also, \ $s'$ \ and \ $t'$ \ have no variables in common.

Let \ $r$ \ be an implicant of \ $p = s'\vee t'$. \ Then \ $r$ \ is an implicant either of \ $s'$ \ or of \ $t'$; \   
for if not, then values can be assigned to variables so that \ $r$ \ is \ $1$ \ while both \ $s'$ \ and \ $t'$ \ are \ $0$, 
\ whence \ $p$ \ will also be \ $0$, \ contradicting the hypothesis that \ $r$ \ is an implicant of \ $p$. \ But if \ $r$ \ 
is a prime implicant of \ $p$, \ then \ $r$ \  cannot be an implicant of both \ $s'$ \ and \ $t'$; \ for then the removal 
from \ $r$ \ of the literals of variables in \ $t$ \ would yield a shorter implicant of \ $s'$ \ and hence of \ $p$. \ 
Therefore the prime implicants of \ $p$ \ are already prime implicants either of \ $s'$ \ or of \ $t'$. \ Claim 2 is proved.

Let \ $p := {\bf u}^\star = st\star = s'\vee t'$ \ with \ ${\bf u}^\star\in F^{\sigma,\star}(k)$. \ Define the binary relation \ $\rho$ \ on \ $\{x_i: i\in k\}$ \ by: \ \ $x_i\rho x_j$ \ iff literals of \ $x_i$ \ and of \ $x_j$ \ appear together in some clause of the complete sum form of \ $p$. \ Claim 1 implies that \ $\rho$ \ is reflexive, and \ $\rho$ \ is symmetric by construction.  Thus the transitive closure, \ $\tau(\rho)$, \ of \ $\rho$ \ is an equivalence relation on \ $\{x_i: i\in k\}$.

We will now prove the following by induction.\vspace{.5 em}

\noindent\underbar{Claim 3}: \ If \ $p =st\star\in F^{\sigma,\star}(k)$ \ for \ $k\ge2$ \ then \ $\tau(\rho)$ \ has exactly 
two equivalence classes.\vspace{.5 em}

Basis: \ \ $p = x_0x_1$. \ The complete sum of this \ $p$ \ is \ $x_0'\vee x_1'$. \ Clearly \ $\{\{x_0\},\{x_1\}\}$ \ 
is the family of equivalence classes of \ $\tau(\rho)$.

For the inductive step, let \ $p = s'\vee t'$, \ and suppose that Claim 3 holds for \ $s$ \ and \ $t$. \ Claim 2 implies that  \ $\tau(\rho)$ \ does not relate variables in \ $s'$ \ with variables in \ $t'$. \ So it remains only to show that all of the variables in \ $s'$, \ say, \ are related to each other by \ $\tau(\rho)$. This is immediate if \ $s$ \ is a single literal. So we may take it that \ $s = u'\vee v'$, \ for formal products \ $u$ \ and \ $v$ \ interpreted in \ $2_{14}$. \ By DeMorgan's Law, \ $s' = u \wedge v$.

We show that the prime implicants of \ $u\wedge v$ \ are precisely the formulas of the form \ $m\wedge n$, \ where \ $m$ \ 
is a prime implicant of \ $u$, \  and \ $n$ \ is a prime implicant of \ $v$: \ Let \ $q$ \ be an implicant of \ $u\wedge v$.  
\ Then \ $q$ \ is an implicant both of \ $u$ \ and of \ $v$. \ Thus \ $q$ \ must be is an implicant of some prime implicant \ $m$ \ of \ $u$ \ and some prime implicant \ $n$ \ of \ $v$. \ So \ $q$ \ must be an implicant of \ $m\wedge n$. \ This shows that every prime implicant of \ $u\wedge v$ \ must be some \ $m\wedge n$. \ But none of the \ $m\wedge n$ \ can imply another; for, suppose \ $m_0\wedge n_0$ \ were an implicant of \ $m_1\wedge n_1$. \ Then \ $m_0 \wedge n_0$ \ is an implicant of \ $m_1$. \  No variables in \ $n_0$ \ appear in \ $m_1$; \ so we can remove their literals, getting that \ $m_0$ \ implies \ $m_1$. \ As prime implicants of \ $u$, \ they are equal.  Similarly \ $n_0 = n_1$. \ We infer the assertion opening this paragraph.

Now let \ $x$ \ be a variable of \ $u$, \ and \ $y$ \ a variable of \ $v$. \ Claim 1 implies that \ $u$ \ depends on \ $x$. \ So \ $x$ \ must appear in some prime implicant \ $m$ \ of \ $u$. \ Similarly, \ $y$ \ appears in some prime implicant \ $n$ \ of \  $v$. \ Thus both \ $x$ \ and \ $y$ \ appear in \ $m\wedge n$, \ which is a prime implicant of \ $s$ \ by the previous paragraph.  Therefore every variable of \ $u$ \ is related by \ $\rho$ \ to every variable of \ $v$. \ So \ $\tau(\rho)$ \ relates all variables in \ $s' = u\wedge v$. \ Claim 3 follows.\vspace{.5em}

For \ $\{{\bf u},{\bf v}\}\subseteq F^\sigma(k)$, \ it is obvious that \ ${\bf u}={\bf v}\,\Rightarrow\,{\bf u^\star}={\bf v}^\star$. \ The converse is obvious for \ $k=1$. \ This is the basis step of an induction on \ $k$.

For the inductive step, let \ $k\in{\bf N}$ \ and suppose for every \ $j\in\{1,2,\ldots,k\}$ \ that \ 
${\bf u}^\star = {\bf v}^\star\,\Rightarrow\,{\bf u} = {\bf v}$ \ whenever \ $\{{\bf u},{\bf v}\}\in F^\sigma(j)$. \ 
Let \ $\{{\bf u},{\bf v}\}\subseteq F^\sigma(k+1)$, \ and suppose that \ ${\bf u}^\star = {\bf v}^\star$. 

As above, in Boolean language we write \ ${\bf u}^\star = p_{\bf u} = s_{\bf u}t_{\bf u}\star = 
s'_{\bf u}\vee t'_{\bf u}$ \ and \ ${\bf v}^\star = p_{\bf v} = s_{\bf v}t_{\bf v}\star = 
s'_{\bf v}\vee t'_{\bf v}$. \ From \ $p_{\bf u} = p_{\bf v}$ \ we get by Claim 3 that \ $s_{\bf u}$ \ has the same 
variables as \ $s_{\bf v}$ \ and that \ $t_{\bf u}$ \ has the same variables has \ $t_{\bf v}$. \ If on \ $2_{14}$ \ 
it happens both that \ $s_{\bf u} = s_{\bf v}$ \ and that \ $t_{\bf u} = t_{\bf v}$ \ then by the inductive hypothesis 
the corresponding factors of \ {\bf u} \ and \ {\bf v} \ in \ $\langle F^\sigma;\odot\rangle$ \ also are equal, and 
therefore \ ${\bf u} = {\bf v}$ \ as alleged.

Without loss of generality, pretend that there is an assignment of values to the variables in \ $s_{\bf u}$ \ which 
gives \ $s_{\bf u}$ \ the value \ $1$ \ while \ $s_{\bf v}$ \ gets the value \ $0$. \ Then, by Claim 1 there is an 
assignment of values to the variables in \ $t_{\bf u}$ \ which gives \ $t_{\bf u}$ \ the value \ $1$. \ It follows for 
these independent value assignments to the elements in \ $\{x_0,x_1,\ldots,x_k\}$ \ that \ $p_{\bf u}$ \ gets the value 
\ $1'\vee 1' = 11\star = 0$ \ while \ $p_{\bf v}$ \ gets the value \ $0'\vee t'_{\bf v} = 0t_{\bf v}\star = 1$, \ 
contrary to the hypothesis that \ ${\bf u}^\star = {\bf v}^\star$.\hfill$\rule{.5em}{.8em}\vspace{1em}$
    
\noindent{\bf Theorem 3.7.} \ { \sl The concrete groupoid \ $2_j$ \ is completely dissociative if and only if \
$j\in\{2,4,8,11,13,14\}$.  }\vspace{1em}

\noindent{\bf Proof.} \ We write \ ${\cal A}\asymp{\cal B}$ \ iff the groupoid \ ${\cal A}$ \ is either 
anti-isomorphic or isomorphic to \ ${\cal B}$. \  Plainly \ $\asymp$ \ is an equivalence relation on \ 
$\underline{\cal G}(2)$. \ The \ $\asymp$ \ equivalence classes of the eight \ $2_j\in\underline{\cal G}(2)$ \
which are non-semigroups are: \ $\{2_2, 2_4, 2_{11}, 2_{13}\}$, \ $\{2_8, 2_{14}\}$, \ and \ $\{2_{10}, 2_{12}\}$. \

Theorem 3.5 gives us that \ $2_{13}$ \ is completely dissociative, and Theorem 3.6 implies that \ $2_{14}$ \ is 
completely dissociative. In \ $2_{10}$, \ the value of an expression depends only on the value of its final input.  
Thus \ $wx\diamond y\diamond z\,\diamond$ \ and \ $wxy\diamond\diamond z\,\diamond$ \ always produce the same value.  
Therefore \ $2_{10}$ \ fails to be $4$--dissociative, and consequently \ $2_{10}$ \ is not completely dissociative.
\hfill$\rule{.5em}{.8em}$

\setlength{\unitlength}{1em}

\begin{picture}(20,11)

\thicklines

\put(10.5,7.4){\line(1,0){7}}

\put(12,2.6){\line(0,1){6}}

\put(11,7.9){$\diamond$}

\put(13,7.8){\bf 0}

\put(15,7.8){\bf 1}

\put(17,7.8){\bf 2}

\put(11,6){\bf 0}

\put(11,4.3){\bf 1}

\put(11,2.6){\bf 2}

\put(13,6){$0$}

\put(15,6){$1$}

\put(17,6){$0$}

\put(13,4.3){$1$}

\put(13,2.7){$0$}

\put(15,4.3){$1$}

\put(17,4.3){$0$}

\put(15,2.7){$0$}

\put(17,2.7){$2$}

\put(14.6,0){\large${\cal D}$}

\end{picture}

\noindent{\bf Theorem 3.8.}  {\sl There are at least seventeen completely dissociative \ $\langle 3;\diamond\rangle$. }\vspace{1em}

\noindent{\bf Proof.} \ The table above exhibits a groupoid, \ ${\cal D}$, \ which we will prove to be completely 
dissociative. During our argument, we will note table entries which we never use. This indicates that \ ${\cal D}$ \ 
is but one of at least seventeen completely dissociative groupoids \ $3_j\in\underline{\cal G}(3)$. \ Our proof uses 
Theorem 3.2.

We note parenthetically the values of \ $\diamond$ \ to which our argument resorts. Let \ $T := \{0,1\}$. \ 
Since \ $0$ \ is idempotent, it is yieldable. (This uses \ $00\diamond = 0$.) \ Since \ $\diamond$ \ is commutative, 
Left Separation is equivalent to Right Separation. Let \ $L_{x,y} := R_{x,y} := \{0\}$ \ for all \ 
$\langle x,y\rangle \in 2^2$, \ (using \ $00\diamond  = 0$ \ and \ $01\diamond = 1$.) 

For Split Separation, note that \ $1^p0^q\vec{g}_{p+q}\in 3^\omega$ \ yields \ $1$ \ when \ $\{p,q\}\subseteq{\bf N}$, \ 
since \ $\{0,1\}$ \ forms a semilattice, (never using the value of \ $01\diamond$.) Similarly \ $0^p2^q\vec{g}_{p+q}$ \ 
yields \ $0$ \ since \ $\{0,2\}$ \ forms a semilattice, (not using the value of \ $20\diamond$.) \ Now suppose \ 
$1\le i<j\le k$. \ Let \ $\vec{g} := 1^i0^{j-i}2^{k-j+1}\vec{g}_{k+1}$. \ If \ ${\bf ab}\odot$ \ is an $i$--split then \ $\vec{g}{\bf ab}\odot^\diamond = 
1^i\vec{g}_i{\bf a}^\diamond 0^{j-i}2^{k-j+1}\vec{g}_{k+1}{\bf b}^\diamond\diamond = 10\diamond = 1$. (This uses cited 
facts and that \ $1$ \ is idempotent.) If \ $\bf{a'b'\odot}$ \ is a $j$--split, then 
$\vec{g}{\bf a'b'\odot}^\diamond = 1^i0^{j-i}\vec{g}_j{\bf a'}^\diamond 2^{k-j+1}\vec{g}_{k+1}{\bf b'}^\diamond\diamond = 12\diamond = 0$, \ (using cited facts, that \ $2$ \ is idempotent, and that \ $12\diamond = 0$.) \ So \ ${\cal D}$ \ 
is completely dissociative.

The values of \ $20\diamond$ \ and \ $21\diamond$ \ were never used in the argument above. So, we can change \ ${\cal D}$ \ 
to make eight other completely dissociative groupoids with \ $\{0\}\not=\{20\diamond,21\diamond\}$.  Since \ ${\cal D}$ \ is 
abelian, we could instead have used \ $\vec{g} = 2^i0^{j-i}1^{k-j+1}\vec{g}_{k+1}$ \ to show Split Separation -- and 
never have used the values of \ $02\diamond$ \ and \ $12\diamond$ \ of \ ${\cal D}$. \ Thus we can make eight other 
completely dissociative groupoids by changing those values in \ ${\cal D}$.\hfill$\rule{.5em}{.8em}$\vspace{1em}

Most of our proofs may be analyzed in the manner above, and slightly modified to produce additional groupoids are 
completely dissociative.\vspace{2em}

\noindent{\large\bf\S4. \ Primitive groupoids}\vspace{.5em}

By the {\em variety} \ ${\bf V}({\cal G})$ \ generated by a groupoid \ ${\cal G}$ \ we mean the closure of \ $\{{\cal G}\}$ 
\ under homomorphic images, subgroupoids and product groupoids of \ ${\cal G}$. \ We will show later that \ ${\cal G}$ \ 
must be completely dissociative if any groupoid in \ ${\bf V}({\cal G})$ \ is. Thus, of special interest are the completely dissociative groupoids which are not forced to be such because of smaller groupoids.

We say that a finite completely dissociative groupoid \ ${\cal P}$ \ is {\sl primitive} iff no smaller groupoid in \ 
${\bf V}({\cal P})$ \ is completely dissociative. 

Observe that all of the $2$--element completely dissociative groupoids are primitive, since the trivial groupoid is a semigroup. We will establish the primitiveness of many other small completely dissociative groupoids.\vspace{1em}

\noindent{\bf Question.} \ Is there a primitive completely dissociative groupoid \ $n_j$ \ for each integer \ 
$n\ge 2$?\vspace{1em}

To proceed with our study of primitive completely dissociative groupoids, we will need a little material from universal
algebra.  For background, we refer the reader is referred to {\bf [5]}, which is a good beginning text and reference.

Our principal tool will be Birkhoff's Theorem, which first appeared in {\bf [1]} and is also carefully developed in {\bf [5]}.  Before stating it, we should first review some terminology.  Everything will be stated for groupoids, although it naturally generalizes to arbitrary algebras.

By a {\em term} we mean an expression built up from variables using the groupoid operation symbol.  Since we are dealing 
only with small terms, we will use infix notation for them in this section.  Examples of terms: \ $x$, \ $x\bullet y$ \ 
and \ $(x\bullet y)\bullet(y\bullet(x\bullet z))$.  An {\em identity} is an equality between terms that is true for all 
values of the variables. It is customary to use \ $\approx$ \ to show that terms are equal in an identity.  We say that an identity {\em holds} in a groupoid iff it is (always) true there, and that an identity holds in a class of groupoids iff it holds in each member of the class.  Alternatively, we can say that a groupoid {\em satisfies} an identity.  Examples of identities are: the Idempotent Law \ ($x\bullet x\approx x$), \ the Commutative Law \ ($x\bullet y\approx y\bullet x$), \ and the Associative Law \ ($x\bullet(y\bullet z)\approx(x\bullet y)\bullet z$).

From this viewpoint, we see that a groupoid is $3$--dissociative if and only if the Associative Law does not hold in it, 
and that a groupoid is completely dissociative if and only if all of the generalizations of the Associative Law fail to hold 
in that groupoid as well.

A {\em variety} is a class of groupoids that is closed under homomorphic images, subgroupoids and (Cartesian) products of elements in that class. If \ $\Sigma$ \ is a set of identities, then the {\em models} of \ $\Sigma$ \ are precisely the groupoids for which all of the identities in \ $\Sigma$ \ hold.  We can now state Birkhoff's Theorem:\vspace{1em}

\noindent{\bf Theorem 4.1.} \ {\sl A class of groupoids is a variety if and only if it is the class of models of a set
of identities.}\vspace{1em} 

We need a related result, which also is due to Birkhoff.\vspace{1em}

\noindent{\bf Theorem 4.2.} \ {\sl If \ ${\cal G}$ \ is an groupoid, then the variety \ ${\bf V}({\cal G})$ \ generated by \ ${\cal G}$ \ is equal to the class of models of the set of all identities holding in \ ${\cal G}$. }\vspace{1em} 

Thus, to show for a groupoid \ ${\cal H}$ \ that \ ${\cal H}\not\in{\bf V}({\cal G})$, \ it suffices to produce an identity that holds in \ ${\cal G}$ \ but does not hold in \ ${\cal H}$. \ So, whenever \ ${\cal G}$ \ fails to be completely dissociative, some generalized associative law is an identity of \ ${\cal G}$. \ By Theorem 4.2, such an identity holds in every groupoid in \ ${\bf V}(G)$. \ So \ ${\cal G}$ \ completely dissociative if  \ ${\bf V}(G)$ \ contains a primitive groupoid.

We are now ready to study primitive completely dissociative groupoids.  As already noted, the \ $6$ \ nonisomorphic completely dissociative groupoids in \ $\underline{\cal G}(2)$ \ are primitive.  What about the groupoid \ ${\cal D}$ \ treated in Theorem 3.8?\vspace{1em}

\noindent{\bf Theorem 4.3.} \ {\sl \ ${\cal D}$ \ is a primitive completely dissociative groupoid. }\vspace{1em} 

\noindent{\bf Proof.} \ We have by 3.8 that \ ${\cal D}$ \ is completely dissociative. Observe that \ ${\cal D}$ \ 
satisfies the Idempotent and Commutative laws. Thus every groupoid in \ ${\bf V}({\cal D})$ \ satisfies them too.
But the only $2$--element groupoids where these laws hold are the semigroups \ $2_{7}$ \ and \ $2_1$. \ This shows 
that there are no completely dissociative groupoids in \ ${\bf V}({\cal D})$ \ which are smaller than \ ${\cal D}$, 
\ and so $\cal D$ is primitive.\hfill$\rule{.5em}{.8em}$\vspace{1em}

We conjecture that the other \ $16$ \ groupoids that were proved completely dissociative in Theorem 3.8 are primitive 
as well.

As another example of our techniques, we will prove that the groupoid \ $\cal B$ \ of Theorem 3.3 is primitive. To  
this end, we investigate the $3$--element groupoids that are commutative and idempotent. The groupoid \ $\cal B$ \ 
satisfies these laws. So all of the groupoids in the variety \ ${\bf V}({\cal B})$ \ also satisfy them.

So, consider the groupoids \ $3_t$ \ that conform to the binary operation table(s) \ $CI3_\alpha$, \ below, with \ 
$\langle a,b,c\rangle\in 3^3$ \ where \ $3 := \{0,1,2\}$. 

\setlength{\unitlength}{1em}

\begin{picture}(20,11)

\thicklines 

\put(10.5,7.4){\line(1,0){7}} 

\put(12,2.6){\line(0,1){6}} 

\put(11,7.9){$\star$}

\put(13,7.8){\bf 0} 

\put(15,7.8){\bf 1} 

\put(17,7.8){\bf 2} 

\put(11,6){\bf 0} 

\put(11,4.3){\bf 1} 

\put(11,2.6){\bf 2} 

\put(13,6){$0$} 

\put(15,6){$a$} 

\put(17,6){$b$} 

\put(13,4.3){$a$} 

\put(13,2.7){$b$} 

\put(15,4.3){$1$} 

\put(17,4.3){$c$} 

\put(15,2.7){$c$} 

\put(17,2.7){$2$}

\put(14.6,0){\large$CI3_{\alpha}$} 
 
\end{picture}\vspace{1.5em}

$CI3_{\alpha}$ \ is our acronym for ``Commutative Idempotent $3$--element groupoid number \ $\alpha$''. \ The index \ 
$\alpha$ \ codes the values of \ $a$, \ $b$ \ and \ $c$ \ according to the following scheme: \ $\alpha = 9a + 3b + c$. 
\ Thus \ $\alpha$ \ ranges from \ $0$ \ to \ $26$.

Many of the \ $CI3_{\alpha}$ \ are isomorphic to each other under permutations of the set \ $3$. \ As one would expect, 
there are $6$--element isomorphism classes, and a few smaller ones. Since proving all the isomorphisms would be tedious,
we merely present the isomorphism classes here:

\begin{itemize}

\item[i)] $CI3_{0} \cong CI3_{13} \cong CI3_{26}$

\item[ii)] $CI3_{1} \cong CI3_{2} \cong CI3_{8}
\cong CI3_{10} \cong CI3_{16} \cong CI3_{17}$

\item[iii)] $CI3_{3} \cong CI3_{12} \cong CI3_{18}
\cong CI3_{22} \cong CI3_{23} \cong CI3_{24}$

\item[iv)] $CI3_{4} \cong CI3_{6} \cong CI3_{9}
\cong CI3_{14} \cong CI3_{20} \cong CI3_{25}$

\item[v)] $CI3_{5} \cong CI3_{15} \cong CI3_{19}$

\item[vi)] $CI3_{7} \cong CI3_{11}$

\item[vii)] $CI3_{21}$
\end{itemize}

The groupoids \ $\cal D$ \ and \ $\cal E$ \ from \S3 are in this list, as are three new completely dissociative groupoids.  
We will examine each isomorphism class briefly, giving tables for one groupoid in each.  We will see that none of the completely dissociative \ $CI3_\alpha$ \ are elements in \ ${\bf V}({\cal B})$, \ since the identity \ ${\bf \beta}$ \ fails 
in each of them, where \ ${\bf \beta}$ \ is: \[    ((x\star y)\star z)\star z\approx((x\star y)\star(x\star z))\star(x\star z)\]

The reader is asked to verify that $\bf{\beta}$ holds in \ ${\cal B}$, \ and hence in \ ${\bf V}({\cal B})$. \ The key to
doing this easily is to note that \ $(x\star u)\star u$ \ is always equal to \ $x$ \ in \ ${\cal B}$, \ except when \ 
$x = 0$ \ and \ $u\neq 0$.

Where it matters, we will indicate how the identity \ ${\bf\beta}$ \ fails.\vspace{1em}

\noindent{i) \ $CI3_{0}\cong CI3_{13}\cong CI3_{26}$}

\setlength{\unitlength}{1em}

\begin{picture}(20,11)

\thicklines 

\put(10.5,7.4){\line(1,0){7}} 

\put(12,2.6){\line(0,1){6}} 

\put(11,7.9){$\star$}

\put(13,7.8){\bf 0} 

\put(15,7.8){\bf 1} 

\put(17,7.8){\bf 2} 

\put(11,6){\bf 0} 

\put(11,4.3){\bf 1} 

\put(11,2.6){\bf 2} 

\put(13,6){$0$} 

\put(15,6){$0$} 

\put(17,6){$0$} 

\put(13,4.3){$0$} 

\put(13,2.7){$0$} 

\put(15,4.3){$1$} 

\put(17,4.3){$0$} 

\put(15,2.7){$0$} 

\put(17,2.7){$2$}

\put(14.6,0){\large$CI3_0$} 
 
\end{picture}\vspace{1.5em}

$CI3_0$ \ is a semigroup. In fact, \ $(x\star y)\star z = 0 = x\star(y\star z)$ \ for every \ $\langle x,y,z\rangle\in  
3^3\setminus\{\langle1,1,1\rangle,\langle2,2,2\rangle\}$. \ If \ $\langle x,y,z\rangle = \langle1,1,1\rangle$ \ then \ 
$(x\star y)\star z = 1 = x\star(y\star z)$, \ and if \ $\langle x,y,z\rangle = \langle 2,2,2\rangle$ \ then \  
$(x\star y)\star z = 2 = x\star(y\star z)$.\vspace{1em}

\noindent{ii) \ $CI3_{1} \cong CI3_{2} \cong CI3_{8}
\cong CI3_{10} \cong CI3_{16} \cong CI3_{17}$}

\setlength{\unitlength}{1em}

\begin{picture}(20,11)

\thicklines 

\put(10.5,7.4){\line(1,0){7}} 

\put(12,2.6){\line(0,1){6}} 

\put(11,7.9){$\star$}

\put(13,7.8){\bf 0} 

\put(15,7.8){\bf 1} 

\put(17,7.8){\bf 2} 

\put(11,6){\bf 0} 

\put(11,4.3){\bf 1} 

\put(11,2.6){\bf 2} 

\put(13,6){$0$} 

\put(15,6){$0$} 

\put(17,6){$0$} 

\put(13,4.3){$0$} 

\put(13,2.7){$0$} 

\put(15,4.3){$1$} 

\put(17,4.3){$1$} 

\put(15,2.7){$1$} 

\put(17,2.7){$2$}

\put(14.6,0){\large$CI3_1$} 
 
\end{picture}\vspace{1.5em}

The groupoid \ $CI3_1$ \ is a semigroup. It is isomorphic to the $3$--element chain under the standard meet operation.\vspace{1em}

\noindent{iii) \ $CI3_{3}\cong CI3_{12}\cong CI3_{18}\cong CI3_{22}\cong CI3_{23}\cong CI3_{24}$}

\setlength{\unitlength}{1em}

\begin{picture}(20,11)

\thicklines 

\put(10.5,7.4){\line(1,0){7}} 

\put(12,2.6){\line(0,1){6}} 

\put(11,7.9){$\star$}

\put(13,7.8){\bf 0} 

\put(15,7.8){\bf 1} 

\put(17,7.8){\bf 2} 

\put(11,6){\bf 0} 

\put(11,4.3){\bf 1} 

\put(11,2.6){\bf 2} 

\put(13,6){$0$} 

\put(15,6){$0$} 

\put(17,6){$1$} 

\put(13,4.3){$0$} 

\put(13,2.7){$1$} 

\put(15,4.3){$1$} 

\put(17,4.3){$0$} 

\put(15,2.7){$0$} 

\put(17,2.7){$2$}

\put(14.6,0){\large$CI3_3$} 
 
\end{picture}\vspace{1.5em}

We verify that the groupoid \ $CI3_3$ \ is completely dissociative. We use Theorem 3.2, with \ $T = 3 = \{0,1,2\}$. \ 
Since \ $2$ \ is idempotent, and its row and column contains every element in \ $3$. \ So we may always set \ $L_{x,y} = R_{x,y} = \{2\}$, \ giving Left and Right Separation.

For Split Separation, let \ $\vec{g} = 0^{i}1^{j-i}2^{k-j+1}\vec{g}_{k+1}\in 3^\omega$. \ If \ ${\bf ab}\odot\in F^\sigma(k+1)$  \ is an $i$--split, then \ $\vec{g}{\bf a}^\star = 0$ \ since \ $0$ \ is idempotent, while \ $\vec{g}_i{\bf b}^\star\in\{0,1\}$, \ since the set \ $\{0,1\}$ \ is absorptive. So \ $\vec{g}{\bf ab}\odot^\star\in\{00\star,01\star\}=\{0\}$. \ But if \ ${\bf a'b'}\odot\in F^\sigma(k+1)$ \ is a $j$--split, then \ $\vec{g}{\bf a'}^\star=0$ \ since \ $\{0,1\}$ \ forms a semilattice, while \ $\vec{g}_j{\bf b'}^\star = 2$ \ because \ $2$ \ 
is idempotent. Thus \ $\vec{g}{\bf a'b'}^\star = 02\star = 1$. \ Therefore \ $CI3_3$ \ is completely dissociative, as alleged.

Since \ $CI3_3$ \ is completely dissociative, since \ $|CI3_3|=3<4=|{\cal B}|$, \ since it is our intention to show that \ $\cal {B}$ \ is primitive, and since the identity \ ${\bf \beta}$ \ holds in \ ${\cal B}$, \ we must show that \ 
$\bf{\beta}$ \ fails to hold in \ $CI3_3$. \ So let  \ $\langle x,y,z\rangle := \langle0,2,1\rangle$, \ and 
observe that then \ $((x\star y)\star z)\star z = ((0\star 2)\star 1)\star 1 = 1\ne 0 = 
((0\star 2)\star(0\star 1))\star(0\star 1) = ((x\star y)\star(x\star z))\star(x\star z)$, \ as desired.\vspace{1em}

\noindent{iv) \ $CI3_{4}\cong CI3_{6}\cong CI3_{9}\cong CI3_{14}\cong CI3_{20}\cong CI3_{25}$}\vspace{.5em}

Since \ $CI3_{9}=\cal{D}$, \ which was proven in Theorem 3.8 to be completely dissociative, it remains to show 
to show that \ $\bf{\beta}$ \ fails in \ ${\cal D}$. \ This time let $\langle x,y,z\rangle = \langle 1,2,0\rangle$, \ 
and note that then \ $((x\star y)\star z)\star z = ((1\star 2)\star 0)\star 0 = 0 \ne 1 = ((1\star 2)\star(1\star 0)) \star(1\star 0) = ((x\star y)\star(x\star z))\star(x\star z)$.\vspace{1em}

\noindent{v) \ $CI3_{5} \cong CI3_{15} \cong CI3_{19}$}\vspace{.5em}

We claim that the groupoid \ $CI3_5$, \ below, is completely dissociative. To use Theorem 3.2, we let \ $T = 3$. \ 
Let \ $L_{x,y} = R_{x,y} = \{1\}$ \ for all \ $\langle x,y\rangle\in3^2$, \ thus establishing Left and Right 
Separation since \ $1$ \ is an identity element of \ $CI3_5$.

\setlength{\unitlength}{1em}

\begin{picture}(20,11)

\thicklines 

\put(10.5,7.4){\line(1,0){7}} 

\put(12,2.6){\line(0,1){6}} 

\put(11,7.9){$\star$}

\put(13,7.8){\bf 0} 

\put(15,7.8){\bf 1} 

\put(17,7.8){\bf 2} 

\put(11,6){\bf 0} 

\put(11,4.3){\bf 1} 

\put(11,2.6){\bf 2} 

\put(13,6){$0$} 

\put(15,6){$0$} 

\put(17,6){$1$} 

\put(13,4.3){$0$} 

\put(13,2.7){$1$} 

\put(15,4.3){$1$} 

\put(17,4.3){$2$} 

\put(15,2.7){$2$} 

\put(17,2.7){$2$}

\put(14.6,0){\large$CI3_{5}$} 
 
\end{picture}\vspace{1.5em}

To show split separation, let \ $1\le i<j\le k$ \ for \ $k\in\{2,3,\ldots\}$, \ and let \ ${\bf ab}\odot$ \ and \ 
${\bf a'b'}\odot$ \ be an $i$--split and a $j$--split, respectively, with \ 
$\{{\bf ab}\odot,{\bf a'b'}\odot\}\subseteq F^\sigma(k+1)$. \ Let \ $\vec{g} := 01^{j-2}2^{k-j+2}\vec{g}_{k+1}\in 3^\omega$.  
\ Then \ $\vec{g}{\bf ab}\odot^\star = 01^{i-1}\vec{g}_i{\bf a}^\star 1^{j-i-1}2^{k-j+2}\vec{g}_{k+1}{\bf b}^\star\star$. \ 
Now, since \ $\{0,1\}$ \ forms a semilattice, and since \ ${\bf a}\in F^\sigma(i)$, \ we have that \ 
$01^{i-1}\vec{g}_i{\bf a}^\star = 0$. \ Likewise, since \ $\{1,2\}$ \ forms a semilattice, we have that \ 
$1^{j-i-1}2^{k-j+2}\vec{g}_{k+1}{\bf b}^\star = 2$. \ Therefore \ $\vec{g}{\bf ab}\odot^\star = 02\star = 1$. \ 
Similarly, \ $\vec{g}{\bf a'b'}\odot^\star = 01^{j-2}2\vec{g}_j{\bf a'}^\star 2^{k-j+1}\vec{g}_{k+1} = c2\star$, \ 
where \ $c := 01^{j-2}2\vec{g}_j{\bf a'}^\star$. \ If \ $j=2$ \ then \ $c = 02\star = 1$, \ whence \ 
$\vec{g}{\bf a'b'}\odot^\star = 12\star = 2$. \ However, if \ $j\ge 3$ \ then \ $c = 01^{j-2}2\vec{g}_j{\bf ps}\odot^\star$ \ for some \ $\langle{\bf p},{\bf s}\rangle\in F^\sigma(t)\times F^\sigma(j-t)$ \ with \ $t\in\{1,2,\ldots,j-1\}$. \ But then \ 
$\vec{g}{\bf p}^\star = 01^{t-1}\vec{g}_t{\bf p}^\star = 0$ \ and \ $\vec{g}_t{\bf s}^\star = 1^{j-t-1}2\vec{g}_j{\bf s}^\star = 2$, \ whence \ $c = \vec{g}{\bf a'}^\star = 02\star = 1$. \ Again \ $\vec{g}{\bf a'b'}^\star = 2$. \ In summary, \ 
$\vec{g}{\bf ab}\odot^\star = 1\not= 2 = \vec{g}{\bf a'b'}\odot^\star$, \ and Split Separation is confirmed.   

Having just proved \ $CI3_5$ \ to be completely dissociative, we must show \ $CI3_5\not\in{\bf V}({\cal B})$. \ 
To this end we show that the identity \ ${\bf\beta}$ \ fails in \ $CI3_5$: \ \ Let \ 
$\langle x,yz\rangle = \langle2,0,1\rangle$. \ Then \ $((x\star y)\star z)\star z = ((2\star 0)\star 1)\star 1 = 1 \not= 2 = 
((2\star 0)\star (2\star 1))\star (2 \star 1) = ((x\star y)\star (x\star z))\star (x\star z)$.\vspace{1em}

\noindent{vi) \ $CI3_{7} \cong CI3_{11}$}

\setlength{\unitlength}{1em}

\begin{picture}(20,11)

\thicklines 

\put(10.5,7.4){\line(1,0){7}} 

\put(12,2.6){\line(0,1){6}} 

\put(11,7.9){$\star$}

\put(13,7.8){\bf 0} 

\put(15,7.8){\bf 1} 

\put(17,7.8){\bf 2} 

\put(11,6){\bf 0} 

\put(11,4.3){\bf 1} 

\put(11,2.6){\bf 2} 

\put(13,6){$0$} 

\put(15,6){$0$} 

\put(17,6){$2$} 

\put(13,4.3){$0$} 

\put(13,2.7){$2$} 

\put(15,4.3){$1$} 

\put(17,4.3){$1$} 

\put(15,2.7){$1$} 

\put(17,2.7){$2$}

\put(14.6,0){\large$CI3_{7}$} 
 
\end{picture}\vspace{1.5em}

We use Theorem 3.2 to prove that \ $CI3_7$ \ is completely dissociative. Let \ $T = \{0,2\}$. \ The element 
\ $0$ \ is idempotent and hence yieldable. Since \ $02\star = 20\star = 2$, \ we have Left and Right Separation 
by always taking \ $L_{x,y} = R_{x,y} = \{0\}$. \ To show Split Separation, let \ $1\le i<j\le k$, \ let \ $\vec{g} := 
0^{i}1^{j-i}2^{k-j+1}$, \ let \ ${\bf ab}\odot$ \ be an $i$--split, and let \ ${\bf a'b'}\odot$ \ be a $j$--split, \ 
where \ $\{{\bf ab}\odot,{\bf a'b'}\odot\}\subseteq F^\sigma(k+1)$. \ Then \ $\vec{g}{\bf ab\odot}^\star = 
0^{i}\vec{g}_i{\bf a}^\star 1^{j-i}2^{k-j+1}\vec{g}_{k+1}{\bf b}^\star\star = 01\star = 0$, \ since \ $\{1,2\}$ \ 
forms a semilattice. Likewise, \ 
$\vec{g}{\bf a'b'\odot}^\star = 0^i1^{j-i}\vec{g}_j{\bf a}^\star 2^{k-j+1}\vec{g}_{k+1}{\bf b}^\star\star =
02\star = 2$, \ since \ $\{0,1\}$ \ forms a semilattice. Therefore \ $\vec{g}{\bf ab}\odot^\star \not= \vec{g}{\bf a'b'}^\star$, \ and Split Separation is established. So \ $CI3_7$ \ is completely dissociative.     

To see that \ $\bf{\beta}$ \ fails in \ $CI3_7$, \ let \ $\langle x,y,z\rangle := \langle 2,1,0\rangle$, \ and note that 
then \ $((x\star y)\star z)\star z = ((2\star 1)\star 0)\star 0 = 0 \not= ((2\star 1)\star(2\star 0))\star(2\star 0) = 
((x\star y)\star(x\star z))\star(x\star z)$.\vspace{1em}

\noindent{\bf Theorem 4.4.} \ {\sl Groupoids which are isomorphic to \ $CI3_{3}$, \ to \ $CI3_{4}$, \ to \ $CI3_{5}$, \ 
or to \ $CI3_{7}$, \ are primitive completely dissociative groupoids.}\vspace{1em}

\noindent{\bf Proof.} \ The argument is identical to that in Theorem 6.3.\hfill$\rule{.5em}{.8em}$\vspace{1em}

\noindent{\bf Corollary 4.5.} \ {\sl \ \ $\cal B$ \ is a primitive completely dissociative groupoid.}\vspace{1em} 

\noindent{\bf Proof.} \ The groupoid \ $\cal B$, \ of Theorem 4.3, satisfies the Idempotent and Commutative laws.
Thus every groupoid in \ ${\bf V}({\cal B})$ \ satisfies them too. But the only $2$--element groupoids where these 
laws hold are the semilattices, \ $2_{7}$ \ and \ $2_1$, \ both of which are semigroups.

The idempotent commutative $3$--element groupoids were studied above. The isomorphism classes of those which are 
completely dissociative -- specifically, those of \ $CI3_{3}$, \ of \ $CI3_{4}$, \ of \ $CI3_{5}$, \ and of \ $CI3_{7}$,
have no elements in common with \ ${\bf V}({\cal B})$, \  since the identity \ $\bf{\beta}$ \ does not hold in them 
but does hold in \ ${\cal B}$.  

This shows that there are no completely dissociative groupoids in \ ${\bf V}({\cal B})$ \ that are smaller than 
${\cal B}$. \ So $\cal B$ is primitive.\hfill$\rule{.5em}{.8em}$\vspace{2em}

\noindent{\large\bf\S5. \ Size sequences}\vspace{.5em}

For \ ${\cal G} := \langle G;\diamond\rangle$ \ an arbitrary groupoid, by Definition {\bf[1}:3.5{\bf]} \ the expression \ 
${\bf SaT}({\cal G})$ \ denotes the integer sequence \ $\langle\,|F^\sigma(k)/{\cal G}|\,\rangle_{k=2}^\infty$. \ This
sequence is called the {\sl subassociativity type} of \ ${\cal G}$.

We extend Definition {\bf[1}:3.7{\bf]}: \ The expression \ $\nu_{{\cal G},k}(i)$ \ denotes the number of \ 
$\approx_{\cal G}$ \ equivalence classes \ $[{\bf u}]_{\cal G}\in F^\sigma(k)/{\cal G}$ \ with \ 
$|[{\bf u}]_{\cal G}|=i$. \ The sequence \ $\langle\,\langle\nu_{{\cal G},k}(i),i\rangle\,\rangle_{i=1}^\infty$ \ is 
called the {\sl $k$--sizing} of \ ${\cal G}$.

In \ {\bf[3]} \ a \ $k$--sizing of \ ${\cal G}$ \ is called a ``size sequence for \ $k$ \ of \ ${\cal G}$".\vspace{.5em}

The \ $k\,$-th \ term of \ ${\bf SaT}({\cal G})$ \ and the $k$--sizing \ of \ ${\cal G}$ \ are related: 
\[ \sum_{i=1}^\infty\nu_{{\cal G},k}(i) = |F^\sigma(k)/{\cal G}|\qquad\mbox{and}\qquad
\sum_{i=1}^\infty i\cdot\nu_{{\cal G},k}(i) = |F^\sigma(k)| = C(k). \]

It may prove fruitful to study the sizings of {\sl subassociative} groupoids; {\it i.e.}, such groupoids as are neither semigroups nor completely dissociative.

$2_{10}\asymp 2_{12}$ \ are the only subassociative elements in \ $\underline{\cal G}(2)$. \
So, up to \ $\asymp$ \ there is exactly one two-element subassociative groupoid.\vspace{.5em}

The groupoid \ ${\cal E} :=\langle 3;\triangle\rangle$, \ below, is $4$--dissociative but not $5$--dissociative, if our 
computer programming is legitimate. \ Thus \ ${\cal E}$ \ is subassociative, and has nontrivial size sequences.\vspace{1em}

\setlength{\unitlength}{1em}

\begin{picture}(20,11)

\thicklines

\put(10.7,9){\line(1,0){7}}

\put(12,4.1){\line(0,1){6}}

\put(10.8,9.4){$\triangle$}

\put(13,9.3){\bf 0}

\put(15,9.3){\bf 1}

\put(17,9.3){\bf 2}

\put(11,7.5){\bf 0}

\put(11,5.8){\bf 1}

\put(11,4.1){\bf 2}

\put(13,7.5){$0$}

\put(15,7.5){$2$}

\put(17,7.5){$1$}

\put(13,5.8){$2$}

\put(13,4.1){$1$}

\put(15,5.8){$1$}

\put(17,5.8){$0$}

\put(15,4.2){$0$}

\put(17,4.2){$2$}

\put(14.7,1.5){\large ${\cal E}$}

\end{picture}

\vspace{1em}

We have learned a little about \ ${\bf SaT}({\cal E})$ \ and the size sequences of \ ${\cal E}$:\vspace{.5em}

$|F^\sigma(3)/{\cal E}| = |F^\sigma(3)| = 2$, \ and so the $3$--sizing of \ ${\cal E}$ \ is \ $\langle\nu_{{\cal
E},3}(1),1\rangle = \langle 2,1\rangle$. \ We list only those terms of a sizing whose first coordinate is greater 
than zero; the other terms we deem ``irrelevant".

$|F^\sigma(4)/{\cal E}| = |F^\sigma(4)| = 5$. \ So the $4$--sizing of \ ${\cal E}$ \ likewise has only one relevant 
term; namely, \ $\langle\nu_{{\cal E},4}(1),1\rangle = \langle 5,1\rangle$.

$|F^\sigma(5)/{\cal E}| = 10$ \ whereas \ $|F^\sigma(5)| = 14$, \ thus repeating our observation that \ ${\cal E}$ \ is 
not $5$--dissociative. The $5$--sizing of \ ${\cal E}$ \ has two relevant terms: \ $\langle 6,1\rangle,\,\langle 4,2\rangle$.

$|F^\sigma(6)/{\cal E}| = 21$ \ whereas \ $|F^\sigma(6)| = 42$. \ The relevant subsequence of the $6$--sizing of \ 
${\cal E}$ \ is \ $\langle 7,1\rangle,\,\langle 7,2\rangle,\,\langle7,3\rangle\,$.\vspace{1em}

We remarked that \ $\Delta$ \ is commutative. Furthermore, every element in \ ${\cal E}$ \ is idempotent. \ It is easy to 
verify that the automorphism group of \ ${\cal E}$ \ is \ Sym$(3)$.\vspace{2em}


\noindent{\large\bf\S6. \ Some \ $\phi:n^k\rightarrow n$ \ are unrepresentable as any \ ${\bf u}^{\vec{\beta}}$. }\vspace{.5em}

The simplest situation, where \ $n=2$ \ and \ $k=3$, \ is the more demanding.\vspace{1em}

\noindent{\bf Lemma 6.1.} \ {\sl There exists a \ $3$--ary operation \ $\phi:2^3\rightarrow 2$ \ such that for no ordered pair \ $\vec{\beta} := \langle\beta_0,\beta_1\rangle$ \ of binary operations \ $\beta_i:2^2\rightarrow 2$ \ does it happen either that \ $x_0x_1x_2\phi = x_0x_1\beta_0x_2\beta_1$ \ or that \ $x_0x_1x_2\phi = x_0x_1x_2\beta_0\beta_1$. }\vspace{1em}

\noindent{\bf Proof.} \ Define \ $\phi:2^3\rightarrow 2$ \ by \ \ $000\phi = 010\phi = 011\phi = 110\phi = 111\phi = 0$, \ \ and \ \ $001\phi = 100\phi = 101\phi = 1$.

The argument consists of four main cases. Two of the cases show that \ $x_0x_1x_2\phi \not= x_0x_1\beta_0x_2\beta_1$ \ while the other two show that \ $x_0x_1x_2\phi \not= x_0x_1x_2\beta_0\beta_1$.

We detail only one case; it will suffice to reveal the nature of our argument.

\underline{Case}: \ \ $00\beta_0 := 1$ \ and \ ${\bf v}^{\langle\beta_0,\beta_1\rangle} := x_0x_1\beta_0x_2\beta_1$.

We show that there is no pair \ $\langle\beta_0,\beta_1\rangle$ \ of binary operations on \ $2$ \ for which \ $\phi = 
{\bf v}^{\langle\beta_0,\beta_1\rangle}$. \ This involves our proceeding step by step through the construction, of the functions \ $\beta_0$ \ and \ $\beta_1$, \ which is mandated by the \ $\phi$ \ specified above and the initial condition \ $00\beta_0 := 1$, \ until we ram into a wall.

From \ $10\beta_1 =: 00\beta_00\beta_1 = 000\phi := 0$, \ we infer that \ $10\beta_1 = 0$. \ Also, \ $11\beta_1 = 00\beta_01\beta_1 = 001\phi := 1$, \ and so \ $11\beta_1 = 1$.

$01\beta_00\beta_1 = 010\phi := 0$ \ provides two possibilities: \ $01\beta_0 = 0$ \ or \ $01\beta_0 = 1$. \ If \ $01\beta_0=1$ \ then \ $11\beta_1 = 01\beta_01\beta_1 = 011\phi := 0$, \ contrary to our prior observation that \ $11\beta_1 = 1$. \ Therefore \ $01\beta_0 = 0$.

Next, \ $01\beta_1 = 01\beta_01\beta_1 = 011\phi := 0$ \ whence \ $01\beta_1 = 0$. \ By \ $10\beta_00\beta_1 = 100\phi := 1$ \ we are again offered two possibilities: \ $10\beta_0 = 0$ \ or \ $10\beta_0 = 1$. \ But if \ $10\beta_0 = 1$ \ then \ $10\beta_1 = 10\beta_00\beta_1 = 100\phi := 1$, \ contrary to our earlier inference that \ $10\beta_1 = 0$. \ Therefore, \ $10\beta_0 = 0$.

Finally, \ $0 = 01\beta_1 = 10\beta_01\beta_1 = 101\phi := 1$, \ and we hit the wall.

We omit the similar second case, which shows that \ $x_0x_1x_2\phi\not= x_0x_1\beta_0x_2\beta_1$ \ when \ $00\beta_0=0$. \ Likewise \ $x_0x_1x_1\phi = x_0x_1x_2\beta_0\beta_1$ \ is impossible.\hfill$\rule{.5em}{.8em}$\vspace{1em}

We used a case-ridden argument to prove Lemma 6.1 \ because there are twice as many formal $3$--ary products interpreted by some duple of binary operations \ $2^2\rightarrow 2$ \ as there are \ $3$--ary operations on the set \ $2:=\{0,1\}$. \ 
However, when either \ $n\ge 3$ \ or \ $k\ge 4$, \ a straightforward counting argument enables us easily to show that the result established for \ $\langle n,k\rangle = \langle 2,3\rangle$ \ extends to every pair \ $\langle n,k\rangle$ \ of 
integers with \ $n\ge 2$ \ and \ $k\ge 3$.\vspace{1em}

\noindent{\bf Theorem 6.2.} \ {\sl For \ $n\ge 2$ \ and \ $k\ge 3$ \ integers, there exists a \ $k$--ary operation \ $\phi:n^k\rightarrow n$ \ such that \ $\phi\not={\bf u}^{\vec{\beta}}$ \ for every \ ${\bf u}\in F^\sigma(k)$ \ and for 
every \ $(k-1)$--tuple \ $\vec{\beta} := \beta_0\beta_1\ldots\beta_{k-2}$ \ of binary operations \ $\beta_i:n^2\rightarrow n$. }\vspace{1em}

\noindent{\bf Proof.} \ Since Lemma 6.1 \ establishes our claim for the case \ $\langle n,k\rangle = \langle 2,3\rangle$, \ 
we may take it that either \ $n\ge 3$ \ or \ $k\ge 4$.

It is clear from  {\bf[6]}, {\bf[7]}, or {\bf[10]} \ that \ $|F^{\sigma,\vec{\beta}}(k)| = C(k-1)$ \ for every 
\ $(k-1)$--tuple \ $\vec{\beta} := \beta_0\beta_1\ldots\beta_{k-2}$ \ of binary operations \ $\beta_j:n^2\rightarrow n$, \ where \ \[ C(k-1) := \frac{1}{2k-1}{2k-1\choose k}\quad\mbox{and \ $C(n)$ \ is the \ $n$\underline{th} \ Catalan number.}\] Since there are \ $n^{n^2(k-1)}$ \ such \ $\vec{\beta}$, \ it follows that the number \ $\Phi(n,k)$ \ of formal $k$--ary products interpreted by some such \ $\vec{\beta}$ \ is \[\Phi(n,k) = n^{n^2(k-1)}C(k-1) = n^{n^2(k-1)}\frac{(2k-2)!}{(k-1)!k!}\,.\] Thus the ratio \ $R(n,k)$ \ of the number \ $n^{n^k}$ \ of distinct \ $k$--ary operations on \ $n$ \ to the number of distinct interpreted formal $k$--products is \[ R(n,k) = \frac{n^{n^k}}{\Phi(n,k)} = n^{n^k-n^2(k-1)}\frac{(k-1)!k!}{(2k-2)!}\,. \]

Notice that \ $R(n,k)>1$ \ for every pair \ $\langle n,k\rangle$ \ of integers such that either \ $n\ge 3$ \ while \ $k\ge 3$ \ or \ $n\ge 2$ \ while \ $k\ge 4$.\hfill$\rule{.5em}{.8em}$\vspace{1em}

In the light of our proof of Theorem 6.2, it seems surprising that, despite the rapid growth of \ $R(n,k)$ \ as 
\ $n$ \ increases, the conclusion of Theorem 6.2 fails strongly for groupoids on infinite universes.\vspace{1em}

\noindent{\bf Theorem 6.3.} \ {\sl Let \ $G$ \ be an infinite set. Then for every \ $k\ge 3$, \ for every \ $\phi:G^k\rightarrow G$, \ and for every \ ${\bf u}\in F^\sigma(k)$, \ there is a $(k-1)$--tuple \ $\vec{\beta} := \beta_0\beta_1\ldots\beta_{k-2}$ \ of binary operations on \ $G$ \ for which \ $\phi={\bf u}^{\vec{\beta}}$. }\vspace{1em}

\noindent{\bf Proof.} \ Our basis for induction is \ $k=3$. \ Let \ ${\bf u} := x_0x_1x_2\bullet\bullet$ \ and \ 
${\bf v} := x_0x_1\bullet x_2\bullet$. \ Of course then \ $F^\sigma(3) = \{{\bf u},{\bf v}\}$. \ Pick \ $\vec{g} 
:= g_0g_1g_2\ldots  \in G^\omega$. \ Let \ $\phi:G^3\rightarrow G$ \ be arbitrary. We use \ ${\bf u}$ \ as our 
paradigm example. There is a bijective binary operation \ $\beta_0:G^2\rightarrow G$. \ Applying \ $\beta_0$ \
to the segment \ $g_1g_2$ \ of \ $\vec{g}$ \ we obtain \ $h := g_1g_2\beta_0\in G$. \ Define \ $\beta_1$ \ so 
that \ $g_0h\beta_1 := \vec{g}\phi$. \ It follows that \ $\vec{g}{\bf u}^{\vec{\beta}} = g_0g_1g_2\beta_0\beta_1 = 
g_0h\beta_1 = \vec{g}\phi$ \ for \ $\vec{\beta} := \beta_0\beta_1$. \ Thus \ ${\bf u}^{\vec{\beta}}=\phi$.

The same trick produces an ordered pair \ $\vec{\gamma}:=\gamma_0\gamma_1$ \ of binary operations on \ $G$ \ 
such that \ ${\bf v}^{\vec{\gamma}} = \phi$. \ The basis is established.\vspace{.5em}

Inductive Step: \ Pick \ $k\ge 3$. \ Suppose for every \ $j\in\{3,4,\ldots,k\}$, \ for every $j$--ary \ function \
$\phi:G^j\rightarrow G$, \ and for every \ ${\bf z}\in F^\sigma(j)$, \ that \ $\phi = z^{\vec{\eta}}$ \ for some
$(j-1)$--tuple \ $\vec{\eta}$ \ of binary operations on \ $G$. \ Pick any \ $\phi:G^{k+1}\rightarrow G$, \ and any 
\ ${\bf w}\in F^\sigma(k+1)$. \ Then \ ${\bf w} = {\bf ab}\odot$ \ with \ 
$\langle{\bf a},{\bf b}\rangle\in F^\sigma(i)\times F^\sigma(k+1-i)$ \ for some \ $i\in\{1,2,\ldots,k\}$. \ Pick \
$\vec{g} := g_0g_1g_2\ldots\in G^\omega$.

Since \ $G$ \ is infinite, there is an $i$--ary bijection \ $\phi':G^i\rightarrow G$. \ Likewise, there is a \ 
$(k+1-i)$--ary bijection \ $\phi'':G^{k+1-i}\rightarrow G$. \ By the inductive hypothesis, there is an \ $(i-1)$--tuple \ $\vec{\beta'} := \beta_0\beta_1\ldots\beta_{i-2}$ \ of binary operations on \ $G$ \ such that \ $\phi' = 
{\bf a}^{\vec{\beta'}}$. \ Likewise, there is a \ $(k-i)$--tuple \ $\vec{\beta''} := \beta_{i-1}\beta_i\ldots\beta_{k-2}$ \ 
of binary operations on \ $G$ \ such that \ $\phi'' = {\bf b}_i^{\vec{\beta''}}$, \ where \ ${\bf b}_i^{\vec{\beta''}}$ \ 
is defined by \ $\vec{g}\,{\bf b}_i^{\vec{\beta''}} := \vec{g}_i{\bf b}^{\vec{\beta''}}$ \ with \ $\vec{g}_i := g_ig_{i+1}g_{i+2}\ldots$ \ Since both \ $\phi'$ \ and \ $\phi''$ \ are bijections onto \ $G$, \ we can define the binary operation \ $\beta_{k-1}:G^2\rightarrow G$ \ such that \ $\vec{g}\phi = \vec{g}\phi'\,\vec{g}_i\phi''\,\beta_{k-1}$. \ So \
$\vec{g}\phi = \vec{g}{\bf a}^{\vec{\beta'}}\,\vec{g}_i{\bf b}^{\vec{\beta''}}\,\beta_{k-1} = \vec{g}{\bf
a}^{\vec{\beta'}}\,\vec{g}\,{\bf b}_i^{\vec{\beta''}}\,\beta_{k-1} = \vec{g}{\bf ab}_i\bullet^{\vec{\beta}} = 
\vec{g}{\bf ab}\odot^{\vec{\beta}} = \vec{g}{\bf w}^{\vec{\beta}}$, \ where \ $\vec{\beta}$ \ is the concatenation \ $\vec{\beta} := \vec{\beta'}\vec{\beta''}\beta_{k-1}$. \ Since \ $\vec{g}\in G^\omega$ \ is arbitrary, \ $\phi = 
{\bf w}^{\vec{\beta}}$.\hfill$\rule{.5em}{.8em}$\vspace{1em}

Each $k$--ary operation \ $\phi:n^k\rightarrow n$ \ is manifested {\it via} an \ $n^k\times (k+1)$ \ matrix \ 
${\bf M}(\phi) = [y_{i,j}]$ \ with entries in \ $n$, \ where the first \ $k$ \ terms of the \ 
$i$\underline{th} row of \ ${\bf M}(\phi)$ \ is the base--$k$ Hindu-Arabic numeral denoting the integer \ $i-1$, \ 
and where \ $y_{i,1}y_{i,2},\ldots,y_{i,k}\phi = y_{i,k+1}$ \ for each \ $i\in\{1,2,\ldots,n^k\}$.\vspace{1em}

\noindent{\bf Question.} \ \ Is there a nice way to identify those matrices \ ${\bf M}(\phi)$ \ 
such that for a given \ ${\bf u}\in F^\sigma(k)$ \ the equality \ $\phi = {\bf u}^{\vec\beta}$ \ 
is satisfied by some $(k-1)$--tuple \ $\vec\beta := \beta_0\beta_1\ldots\beta_{k-2}$ \ of binary 
operations \ $\beta_i:n^2\rightarrow n\,$?\vspace{1em}

It is reasonable to wonder whether enlarging our tool kit of building-block operations on \ $n$ \ enables the construction of all operations of given arities larger than the arities of permitted building blocks. In this light we ask\vspace{1em}

\noindent{\bf Question.} \ For each \ $r\in\{3,4,5,\ldots\}$, \ is there an \ $n(r)\in{\bf N}$ \ such that, for each pair \ $\langle m,k\rangle$ \ of integers with \ $m\ge n(r)$ \ and \ $k\ge r$, \ there is some $k$--ary operation \ $\phi:m^k\rightarrow m$ \ which it is impossible to ``build" using a natural formal product construction generalizing \ $F^\sigma$ \ by allowing \ $j$--ary operations on \ $m$ \ with \ $j\in\{2,3,\ldots,r-1\}$ \ instead of using only binary operations?\vspace{2em}

\noindent{\large\bf\S7. \ Minimally associative groupoids.}\vspace{.5em}

A groupoid \ $\langle G;\diamond\rangle$ \ is said to be
$k$--{\sl anti-associative} iff for every \ ${\bf u}\not={\bf v}$ \ with \
$\langle{\bf u},{\bf v}\rangle\in F^\sigma(k)\times F^\sigma(k)$ \
and for every \ $\vec{g}\in G^\omega$ \ we have \
$\vec{g}{\bf u}^\diamond\not=\vec{g}{\bf v}^\diamond$.

We deem a finite groupoid  \ $\langle n;\diamond\rangle$ \ to be \
$\langle n,k,m\rangle$--{\sl anti-associative} iff for each pair \
${\bf u}\not={\bf v}$ \ of formal $k$--products we have \
$|\{\vec{g}: \vec{g}{\bf u}^\diamond = \vec{g}{\bf v}^\diamond\}|\le m$.

Thus we see that a groupoid \ $\langle n;\diamond\rangle$ \ is
$k$--anti-associative if and only if \ $\langle n;\diamond\rangle$ \ is \
$\langle n,k,0\rangle$--anti-associative.\vspace{.5em}

Theorem 4.2 in {\bf[2]} states that \ $\langle F^\sigma;\odot\rangle$ \ is $k$--anti-associative 
for each \ $k\ge 3$. \ A similar example of a groupoid that is $k$--anti-associative for all \ $k\ge 3$ \ is 
the free groupoid on one or more generators.  (See {\bf [5]} for definitions and relevant theorems.)  Both of 
these examples are infinite groupoids. 

By the Pigeonhole Principle, if \ $C(k-1)>n$ \ then \ $\langle n;\diamond\rangle$ \ fails to be
$k$--anti-associative, where \ $|F^\sigma(k)| = C(k-1)$ \ is the \ $(k-1)\,$-st \ Catalan number. 
So, no finite groupoid is $k$--anti-associative for all \ $k\ge 3$. \ But, at least for small \ $k$, \ 
there are finite 
$k$--anti-associative groupoids; {\it e.g.}, both \ $2_{10}$ \ and \ $2_{12}$ \ are $3$--anti-associative.  
In {\bf [3]} we will investigate further the\vspace{1em}

\noindent{\bf Question.} \ For each \ $k\geq 3$, \ is there a $k$--anti-associative \ $n_j$?\vspace{1em} 

In any event, for each \ $\langle n,k\rangle$ \ there is a smallest integer \
$M(n,k)$ \ such that there
exists an \ $\langle n,k,M(n,k)\rangle$--anti-associative groupoid \ $\langle n;\alpha\rangle$. \ 
We say that any such groupoid \ $\langle n;\alpha\rangle$ \ is {\sl minimally $k$--associative}.\vspace{1em}

\noindent{\bf Question.} \ \ Characterize the function \ $M:\langle n,k\rangle\mapsto M(n,k)$ \ and the family 
of all minimally$k$--associative groupoids \ $\langle n;\alpha\rangle$.\vspace{2em}

\noindent{\large\bf Acknowledgments.} \ Conversations with Milton Proc\'opio de Borba, Sohan Perdomo Moreira, and Sylvia B.
Silberger assisted our work. We are indebted to Universidade Federal de Santa Catarina for its hospitality, and in particular to its Department Chair, Nereu E. Burin, and to the Head of its PET project, Jos\'e L. R. Pinho, for their having
facilitated our collaboration.\vspace{2.5em}

\noindent{\large\bf References.}\vspace{.5em}

\noindent{\bf [1]} \ \ G. Birkhoff: {\sl On the structure of abstract algebras} Proc. Camb. Philos. Soc. {\bf 31} 
(1935), 433 - 454.\vspace{.5em}

\noindent{\bf [2]} \ \ Milton S. Braitt and Donald Silberger: {\sl Subassociative groupoids}, Quasigroups and Related 
Systems {\bf 14} (2006), 11 - 26.\vspace{.5em} 

\noindent{\bf [3]} \ \ Milton S. Braitt, David Hobby and Donald Silberger: {\sl Anti-associative groupoids}. Preprint available.\vspace{.5em}

\noindent{\bf [4]} \ \ Milton S. Braitt, David Hobby and Donald Silberger: {\sl Size sequences of a groupoid}. In preparation.\vspace{.5em}

\noindent{\bf [5]} \ \ Stanley Burris and H. P. Sankappanavar: {\sl A course in universal algebra}, Springer Verlag, 
1981. (Also freely available online at:

{\nobreak http://www.math.uwaterloo.ca/\verb|~|snburris/htdocs/ualg.html)}\vspace{.5em}

\noindent{\bf [6]} \ \ H. W. Gould: {\sl Research Bibliography of Two Special Sequences}, Combinatorial Research 
Institute, West Virginia University, Morgantown, 1977.\vspace{.5em} 

\noindent{\bf [7]} \ \ I. M. Niven: {\sl Mathematics of Choice or How to Count without Counting}, final chapter. 
New Mathematical Library, Vol. {\bf 15}, 1965. Random House, New York.\vspace{.5em} 

\noindent{\bf [8]} \ \ W. V. Quine:  {\sl A way to simplify truth functions}, American Mathematical Monthly, 
{\bf 62} No. 9 (1955), 627 - 631.\vspace{.5em}

\noindent{\bf [9]} \ \ W. V. Quine:  {\sl Selected logic papers}, Random House, New York, 1966.\vspace{.5em}

\noindent{\bf [10]} \ \ D. M. Silberger: {\sl Occurrences of the integer $(2n-2)!/n!(n-1)!$}, Prace Mat. 
{\bf 13} (1969), 91 - 96.\vspace{2em}

\noindent{\large\bf Addresses.}\vspace{.5em}

Donald Silberger or David Hobby, Department of Mathematics, State University of New York at New Paltz, NY 12561, U.S.A.

Email: \ \ DonaldSilberger@gmail.com \ \ or \ \ hobbyd@newpaltz.edu\vspace{.5em}

Milton S. Braitt, Departamento de Matem\'atica, Universidade Federal de Santa Catarina, Cidade Universit\'aria, 
Florian\'opolis, SC 88040-900, Brasil

Email: \ \ MSBraitt@mtm.ufsc.br\vspace{2em}

\noindent{\bf 2000 Mathematics Subject Classification:}

Primary: \ 20N02
  
Secondary: \ 05A99, 08A99, 08B99, 08C10 \vspace{2em}

\noindent{\bf Keywords:} \ groupoid, generalized associative, dissociative, anti-associative 

\end{document}